\documentclass[%
 reprint,
 amsmath,amssymb,
 aps,
 prx,
]{revtex4-2}

\usepackage{graphicx}
\usepackage{dcolumn}
\usepackage{bm}
\usepackage{mathrsfs} 
\usepackage{xcolor,float}

\usepackage{subcaption}

\usepackage{amsmath}
\usepackage{amsthm}
\usepackage{amssymb}
\usepackage{amsfonts}
\usepackage{mathtools}
\usepackage{enumerate}
\usepackage{comment}
\allowdisplaybreaks

\usepackage{acronym}
\usepackage{latexsym}
\usepackage{paralist}
\usepackage{wasysym}
\usepackage{xspace}

\usepackage[font=small,labelfont=bf]{caption}
\usepackage{tikz}
\usetikzlibrary{calc,patterns}

\usepackage{hyperref}

\numberwithin{equation}{section}  
\usepackage[sort&compress,capitalize,nameinlink]{cleveref}
\crefname{app}{Appendix}{Appendices}

\crefrangeformat{equation}{\upshape(#3#1#4)\textendash(#5#2#6)}

\usepackage{enumitem}

\newcommand{\debug}[1]{{\color{black}#1}}

\begin{document}
\theoremstyle{plain}
\newtheorem{theorem}{Theorem}
\newtheorem{corollary}[theorem]{Corollary}
\newtheorem*{corollary*}{Corollary}
\newtheorem{lemma}[theorem]{Lemma}
\newtheorem{proposition}[theorem]{Proposition}
\newtheorem{conjecture}[theorem]{Conjecture}

\theoremstyle{definition}
\newtheorem{definition}[theorem]{Definition}
\newtheorem*{definition*}{Definition}
\newtheorem{assumption}[theorem]{Assumption}
\newtheorem{hypothesis}{Hypothesis}
\newtheorem*{hypothesis*}{Hypothesis}
\newtheorem{notation}{Notation}

\theoremstyle{remark}
\newtheorem{remark}[theorem]{Remark}
\newtheorem*{remark*}{Remark}
\newtheorem*{notation*}{Notational remark}
\newtheorem{example}[theorem]{Example}


\numberwithin{theorem}{section}


\providecommand\given{} 

\DeclarePairedDelimiter{\braces}{\{}{\}}
\DeclarePairedDelimiter{\bracks}{[}{]}
\DeclarePairedDelimiter{\parens}{(}{)}

\DeclarePairedDelimiter{\abs}{\lvert}{\rvert}
\DeclarePairedDelimiter{\norm}{\lVert}{\rVert}
\newcommand{\dnorm}[1]{\norm{#1}_{\ast}}

\DeclarePairedDelimiter{\ceil}{\lceil}{\rceil}
\DeclarePairedDelimiter{\floor}{\lfloor}{\rfloor}
\DeclarePairedDelimiter{\clip}{[}{]}
\DeclarePairedDelimiter{\negpart}{[}{]_{-}}
\DeclarePairedDelimiter{\pospart}{[}{]_{+}}

\DeclarePairedDelimiter{\bra}{\langle}{\rvert}
\DeclarePairedDelimiter{\ket}{\lvert}{\rangle}
\DeclarePairedDelimiterX{\braket}[2]{\langle}{\rangle}{#1,#2}

\DeclarePairedDelimiterX{\inner}[2]{\langle}{\rangle}{#1,#2}
\DeclarePairedDelimiterX{\setdef}[2]{\{}{\}}{#1:#2}

\DeclarePairedDelimiterXPP{\probof}[1]{\Prob}{(}{)}{}{%
\renewcommand\given{\nonscript\,\delimsize\vert\nonscript\,\mathopen{}}
#1}

\DeclarePairedDelimiterXPP{\exof}[1]{\Expect}{[}{]}{}{%
\renewcommand\given{\nonscript\,\delimsize\vert\nonscript\,\mathopen{}}
#1}
\newcommand{\FC}[1]{\textcolor{red}{\small #1\\}}


\newcommand{\diff}{\ \textup{\debug d}}

\newcommand{\naturals}{\mathbb{\debug N}}
\newcommand{\reals}{\mathbb{\debug R}}
\newcommand{\ident}{\debug I}
\newcommand{\ortho}{\debug U}
\newcommand{\jordan}{\debug J}
\newcommand{\genmatr}{\debug M}
\newcommand{\simplex}{\debug \Sigma}
\newcommand{\run}{\debug k}
\newcommand{\h}{\boldsymbol{h}}
\newcommand{\w}{\boldsymbol{w}}
\newcommand{\x}{\boldsymbol{x}}
\newcommand{\y}{\boldsymbol{y}}

\newcommand{\Expect}{\mathbf{\debug E}}
\newcommand{\Prob}{\mathbf{\debug P}}
\newcommand{\Prw}{\mathsf{\debug W}}
\newcommand{\initial}{\debug \mu}
\newcommand{\transit}{\debug \pi}
\newcommand{\transitvec}{\boldsymbol{\transit}}
\newcommand{\transitvecalt}{\transitvec'}
\newcommand{\transitmatrix}{\debug \Pi}
\newcommand{\transitmatrixalt}{\transitmatrix'}
\newcommand{\transits}{\debug {\mathcal{T}}}
\newcommand{\stationary}{\debug \rho}
\newcommand{\laplace}{\debug L}
\newcommand{\eigen}{\debug \lambda}
\newcommand{\limittransitmatrix}{\overline{\debug \Pi}}
\newcommand{\proba}{\debug p}


\newcommand{\horizon}{\debug T}
\newcommand{\per}{\debug t}

\newcommand{\graph}{\mathcal{\debug G}}
\newcommand{\graphalt}{\graph'}
\newcommand{\vertices}{\mathcal{\debug V}}
\newcommand{\edge}{\debug e}
\newcommand{\edges}{\mathcal{\debug E}}
\newcommand{\edgesalt}{\edges'}
\newcommand{\edgessub}{\mathcal{\debug C}}
\newcommand{\cycle}{\debug \gamma}
\newcommand{\cyclealt}{\debug \xi}
\newcommand{\cycles}{\mathcal{\debug S}}
\newcommand{\cyclealts}{\mathcal{\debug X}}
\newcommand{\unicycle}{\debug \kappa}
\newcommand{\unicycles}{\mathcal{\debug K}}
\newcommand{\outdegr}{\debug \delta}
\newcommand{\tree}{\debug \tau}
\newcommand{\trees}{\mathcal{\debug T}}
\newcommand{\neighbors}{\mathcal{\debug N}}
\newcommand{\neighborsout}{\mathcal{\debug N}^{+}}
\newcommand{\neighborsin}{\mathcal{\debug N}^{-}}
\newcommand{\completegraph}{\debug K}
\newcommand{\cyclegraph}{\debug C}
\newcommand{\ringstargraph}{\debug R}
\newcommand{\arborescence}{\debug \alpha}
\newcommand{\arborescences}{\mathcal{\debug A}}

\newcommand{\game}{\debug \Gamma}

\newcommand{\play}{\debug i}
\newcommand{\playalt}{\debug j}
\newcommand{\nplayers}{\debug n}
\newcommand{\players}{\bracks{\nplayers}}

\newcommand{\pure}{\alpha}
\newcommand{\purealt}{\beta}
\newcommand{\nPures}{\debug S}
\newcommand{\pures}{\mathcal{\nPures}}
\newcommand{\peq}{\pure^{\ast}}

\newcommand{\act}{\debug s}
\newcommand{\actalt}{\debug \act'}
\newcommand{\actprof}{\boldsymbol{\act}}
\newcommand{\actprofalt}{\actprof'}
\newcommand{\actions}{\debug \Sigma}
\newcommand{\acts}{\actions}
\newcommand{\intacts}{\intfeas}
\newcommand{\unactions}{\debug {\widetilde{\actions}}}

\newcommand{\mixed}{\sigma}
\newcommand{\mixedprof}{\boldsymbol{\mixed}}

\newcommand{\cost}{\debug c}
\newcommand{\costprof}{\boldsymbol{\cost}}
\newcommand{\Cost}{C}
\newcommand{\costs}{\debug {\mathcal{C}}}

\newcommand{\pay}{\debug u}
\newcommand{\payv}{v}
\newcommand{\loss}{\ell}

\newcommand{\reg}{R}

\newcommand{\equilibria}{\mathcal{E}}
\newcommand{\correquilibria}{\mathcal{C}}

\newcommand{\cor}[1]{#1^{\ast}}
\newcommand{\eq}[1]{#1^{\ast}}
\newcommand{\optt}[1]{\tilde#1}
\newcommand{\out}[1]{\tilde#1}

\newcommand{\potential}{\debug \Psi}


\newcommand{\argdot}{\,\cdot\,}
\newcommand{\dkl}{D_{\textup{KL}}}
\newcommand{\filter}{\mathcal{F}}
\newcommand{\gen}{\mathcal{L}}
\newcommand{\interval}{\debug I}
\newcommand{\intsimplex}{\simplex^{\!\circ}}
\newcommand{\set}{\mathcal{S}}
\newcommand{\step}{\gamma}
\newcommand{\temp}{\eta}


\newcommand{\depend}{\debug q}

\newcommand{\si}{\sigma} 
\newcommand{\ent}{{\rm ENT} } 
\newcommand{\var}{{\rm Var} } 
\newcommand{\xy}{{x,y} } 
\newcommand{\wt}{\widetilde } 
\newcommand{\tc}{\, |\, } 
\newcommand{\Hd}{H^{low}}
\newcommand{\Hu}{H^{up}}
\newcommand{\whp}{\textbf{whp}}
\newcommand{\ind}{\mathds{1}}
\newcommand{\p}{\mathfrak{p}}
\newcommand{\tx}{\texttt{tx}}

\newcommand{\cA}{\ensuremath{\mathcal A}} 
\newcommand{\cB}{\ensuremath{\mathcal B}} 
\newcommand{\cC}{\ensuremath{\mathcal C}} 
\newcommand{\cD}{\ensuremath{\mathcal D}} 
\newcommand{\cE}{\ensuremath{\mathcal E}} 
\newcommand{\cF}{\ensuremath{\mathcal F}} 
\newcommand{\cG}{\ensuremath{\mathcal G}} 
\newcommand{\cH}{\ensuremath{\mathcal H}} 
\newcommand{\cI}{\ensuremath{\mathcal I}} 
\newcommand{\cJ}{\ensuremath{\mathcal J}} 
\newcommand{\cK}{\ensuremath{\mathcal K}} 
\newcommand{\cL}{\ensuremath{\mathcal L}} 
\newcommand{\cM}{\ensuremath{\mathcal M}} 
\newcommand{\cN}{\ensuremath{\mathcal N}} 
\newcommand{\cO}{\ensuremath{\mathcal O}} 
\newcommand{\cP}{\ensuremath{\mathcal P}} 
\newcommand{\cQ}{\ensuremath{\mathcal Q}} 
\newcommand{\cR}{\ensuremath{\mathcal R}} 
\newcommand{\cS}{\ensuremath{\mathcal S}} 
\newcommand{\cT}{\ensuremath{\mathcal T}} 
\newcommand{\cU}{\ensuremath{\mathcal U}} 
\newcommand{\cV}{\ensuremath{\mathcal V}} 
\newcommand{\cW}{\ensuremath{\mathcal W}} 
\newcommand{\cX}{\ensuremath{\mathcal X}} 
\newcommand{\cY}{\ensuremath{\mathcal Y}} 
\newcommand{\cZ}{\ensuremath{\mathcal Z}} 


\newcommand{\bbA}{{\ensuremath{\mathbb A}} } 
\newcommand{\bbB}{{\ensuremath{\mathbb B}} } 
\newcommand{\bbC}{{\ensuremath{\mathbb C}} } 
\newcommand{\bbD}{{\ensuremath{\mathbb D}} } 
\newcommand{\bbE}{{\ensuremath{\mathbb E}} } 
\newcommand{\bbF}{{\ensuremath{\mathbb F}} } 
\newcommand{\bbG}{{\ensuremath{\mathbb G}} } 
\newcommand{\bbH}{{\ensuremath{\mathbb H}} } 
\newcommand{\bbI}{{\ensuremath{\mathbb I}} } 
\newcommand{\bbJ}{{\ensuremath{\mathbb J}} } 
\newcommand{\bbK}{{\ensuremath{\mathbb K}} } 
\newcommand{\bbL}{{\ensuremath{\mathbb L}} } 
\newcommand{\bbM}{{\ensuremath{\mathbb M}} } 
\newcommand{\bbN}{{\ensuremath{\mathbb N}} } 
\newcommand{\bbO}{{\ensuremath{\mathbb O}} } 
\newcommand{\bbP}{{\ensuremath{\mathbb P}} } 
\newcommand{\bbQ}{{\ensuremath{\mathbb Q}} } 
\newcommand{\bbR}{{\ensuremath{\mathbb R}} } 
\newcommand{\bbS}{{\ensuremath{\mathbb S}} } 
\newcommand{\bbT}{{\ensuremath{\mathbb T}} } 
\newcommand{\bbU}{{\ensuremath{\mathbb U}} } 
\newcommand{\bbV}{{\ensuremath{\mathbb V}} } 
\newcommand{\bbW}{{\ensuremath{\mathbb W}} } 
\newcommand{\bbX}{{\ensuremath{\mathbb X}} } 
\newcommand{\bbY}{{\ensuremath{\mathbb Y}} } 
\newcommand{\bbZ}{{\ensuremath{\mathbb Z}} } 


\newcommand{\E}{\ensuremath{\mathbb{E}}}
\newcommand{\XX}{\ensuremath{\mathcal{X}}}
\newcommand{\C}{\ensuremath{\mathbb{C}}}
\newcommand{\N}{\ensuremath{\mathbb{N}}}
\newcommand{\Z}{\ensuremath{\mathbb{Z}}}
\newcommand{\Q}{\ensuremath{\mathbb{Q}}}
\newcommand{\R}{\ensuremath{\mathbb{R}}}
\renewcommand{\P}{\ensuremath{\mathbb{P}}}
\newcommand{\W}{\ensuremath{\mathbb{W}}}
\newcommand{\mW}{\ensuremath{\widetilde{\mathbb{W}}}}
\newcommand{\pr}{\ensuremath{\right\rangle}}
\renewcommand{\H}{\ensuremath{\mathcal{H}}}
\newcommand{\U}{\ensuremath{\mathcal{U}}}
\newcommand{\K}{\ensuremath{\mathcal{K}}}
\newcommand{\BB}{\ensuremath{\mathcal{B}}}
\newcommand{\FF}{\ensuremath{\mathcal{F}}}
\newcommand{\EE}{\ensuremath{\mathcal{E}}}
\newcommand{\GG}{\ensuremath{\mathcal{G}}}
\newcommand\independent{\protect\mathpalette{\protect\independenT}{\perp}}
\def\independenT#1#2{\mathrel{\rlap{$#1#2$}\mkern2mu{#1#2}}}


%
\def\gap{\mathop{\rm gap}\nolimits}


\preprint{APS/123-QED}



%
\title {Voter model on heterogeneous directed networks}
\thanks{
The work of L.A., R.S.H. and F.C. is supported in part by the Netherlands Organisation for Scientific Research (NWO) through the Gravitation {\sc Networks} grant 024.002.003. The work of F.C. is further supported by the European Union's Horizon 2020 research and innovation programme under the Marie Sk\l odowska-Curie grant agreement no.\ 945045. 
The work of D.G. is supported by the European Union - NextGenerationEU - National Recovery and Resilience Plan (Piano Nazionale di Ripresa e Resilienza, PNRR), projects `SoBigData.it - Strengthening the Italian RI for Social Mining and Big Data Analytics' - Grant IR0000013 (n. 3264, 28/12/2021)
and ``Reconstruction, Resilience and Recovery of Socio-Economic Networks'' RECON-NET EP\_FAIR\_005 - PE0000013 ``FAIR'' - PNRR M4C2 Investment 1.3, financed by the European Union – NextGenerationEU.
This work was performed using the compute resources from the Academic Leiden Interdisciplinary Cluster Environment (ALICE) provided by Leiden University.\hfill
\parbox{0.1\textwidth}
{~~~~\includegraphics[width=0.05\textwidth]{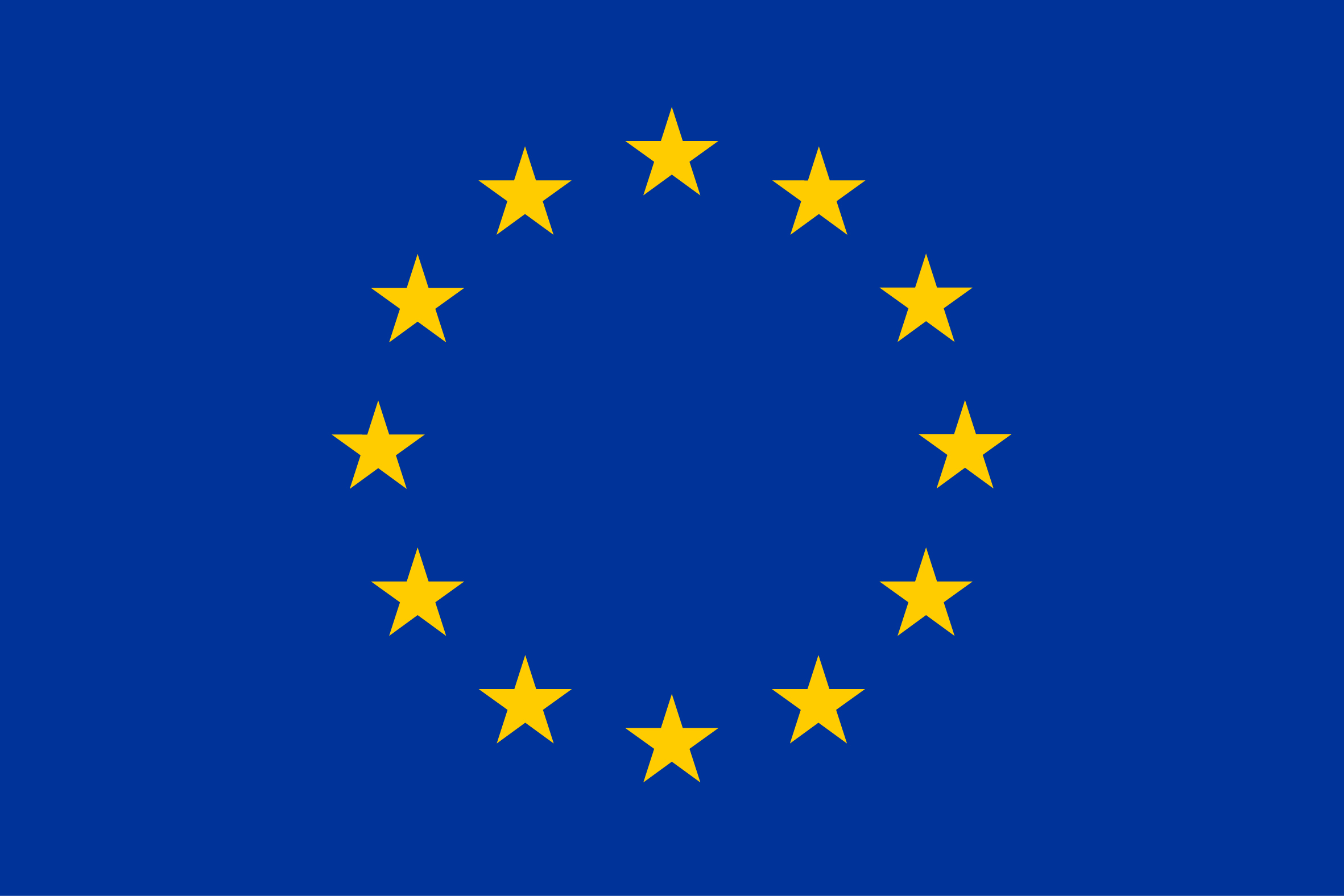}}}%

 \author{Luca Avena}
 \affiliation{Dipartimento di Matematica e Informatica "Ulisse Dini", Università degli Studi di Firenze, Viale Morgagni 67/A, 50134, Firenze, Italy}
 \author{Federico Capannoli}
 \affiliation{Mathematical Institute, Leiden University, Gorlaeus Gebouw, BW-vleugel, Einsteinweg 55, 2333 CC Leiden, The Netherlands}
 \author{Diego Garlaschelli}
 \affiliation{IMT School for Advanced Studies, P.zza San Francesco 19, 55100 Lucca, Italy}
 \affiliation{Lorentz Institute for Theoretical Physics, University of Leiden, Niels Bohrweg 2, 2333 CA Leiden, The Netherlands}
 \author{Rajat Subhra Hazra}
 \affiliation{Mathematical Institute, Leiden University, Gorlaeus Gebouw, BW-vleugel, Einsteinweg 55, 2333 CC Leiden, The Netherlands}




\begin{abstract}

We investigate the consensus dynamics of the voter model on large random graphs with heterogeneous and directed features, focusing in particular on networks with power-law degree distributions. By extending recent results on sparse directed graphs, we derive exact first-order asymptotics for the expected consensus time in directed configuration models with i.i.d. Pareto-distributed in- and out-degrees. For any tail exponent $\alpha>0$, we derive the mean consensus time scaling depending on the network size and a pre-factor that encodes detailed structural properties of the degree sequences. We give an explicit description of the pre-factor in the directed setting. This extends and sharpens previous mean-field predictions from statistical physics, providing the first explicit consensus-time formula in the directed heavy-tailed setting.
Through extensive simulations, we confirm the validity of our predictions across a wide range of heterogeneity regimes, including networks with infinite-variance and infinite-mean degree distribution. We further explore the interplay between network topology and voter dynamics, highlighting how degree fluctuations and maximal degrees shape the consensus landscape. Complementing the asymptotic analysis, we provide numerical evidence for the emergence of Wright–Fisher diffusive behavior in both directed and undirected ensembles under suitable mixing conditions, and demonstrate the breakdown of this approximation in the in the infinite mean regime. 
\end{abstract}

\maketitle       
\section{Introduction}
The voter model is a classical interacting particle system that models how consensus (unanimous opinion) emerges on a network $G= (V, E)$. In the standard version on a finite connected graph each vertex holds one of two opinions and, at exponential random times, a randomly chosen vertex adopts the opinion of a random neighbor.
The voter model was originally studied in \cite{CS73, HL75} and can be described through a continuous-time Markov chain with state space $\{0,1\}^V$ (see Section \ref{sec:VoterDef} for a more precise definition).
If the underlying graph is infinite graph the characterization of the set of invariant measures of the voter model is non-trivial (see for example \cite{Lig85}). Whereas if the underlying graph is finite and strongly connected, that is, there exists a directed path from any two pair of vertices, then
the set of possible invariant measures become trivial and is expressed in terms of  two Dirac measures on the monochromatic configurations in $\{0,1\}^V$ made of all $1$'s and all $0$'s. Thus, on finite connected graphs the voter model trivially achieves \emph{consensus}, that is, one of the two monochromatic  configurations, in finite time. One of the main interest is in understanding how consensus is reached as a function of the size of the graph and how this is related to the geometry of the underlying graph. A key quantity is the \emph{consensus time} $\tau^u_{\rm cons}$, that is, the first time the system reaches agreement. Throughout the whole paper we will denote by $\tau^u_{\rm cons}$ the consensus time for the voter model where the initial distribution of the opinions is given by a product of $\text{Bern}(u)$ random variables, $u\in[0,1]$.

The study of consensus time in the voter model has inspired a rich and interdisciplinary literature across mathematics and physics. Early foundational work emerged in the 1980s and 1990s (see \cite{Lig85, AldPCH, Dur07}), laying the groundwork for rigorous analysis on regular graphs. A canonical example is the complete graph, where the voter model coincides with the classical \emph{Moran model} from population genetics (see \cite{Mor58}). In this setting, the fraction of individuals holding a particular opinion evolves as a \emph{Wright--Fisher diffusion}, and the expected time to consensus grows linearly with the population size $n$. The exact distribution and expected value of the consensus time are well understood in this case, making it a fundamental benchmark for comparison.

There are many variants of the voter model that have been extensively studied in the computer science and physics literature; we present here a brief summary.
The classical voter model studied here involves a node adopting the opinion of a random neighbor, known as the \emph{pull} model.  
Alternatives include the \emph{push} voter model \cite{Dur10}, the \emph{oblivious} model \cite{CDFR18}, and time evolution variants such as synchronous and asynchronous updates.  
Further generalizations introduce randomness (the \emph{noisy} model) \cite{GM95}, opinion thresholds like the \( q \)-voter \cite{DS93, CMP09} and \( k \)-majority models \cite{DeOli92}.  
Biases are also studied, including single-bias majority dynamics \cite{HD19, cruciani2023phase}.  
For a thorough overview, see Rivera’s PhD thesis \cite{Riv18thesis} and related applied works \cite{hassin1999distributed, YPOS10, BCP11}.

Beyond fully connected graphs, exact asymptotics for consensus time have been derived in more structured geometries. For instance, \cite{cox:greven, Cox89} analyzed the voter model on the $d$ dimensional torus of volume $n$, proving that the consensus time scales as $ n^{2/d}$ in $d=1$, $n^{2/d}\log n$ in $d=2$ and $n$ in $d\ge 3$ as the system size tends to infinity.  On the other hand, \cite{CF05} showed that in the random $d$-regular graph ensemble---the consensus time, after proper rescaling, converges in probability to a deterministic constant, $(d-1)/(d-2)$. Furthermore, the density of a given opinion evolves according to the Wright--Fisher diffusion on the consensus time scale, much like in the Moran model (\cite{CCC16}). These developments have motivated extensive investigation into how graph topology, degree heterogeneity, and directedness influence consensus dynamics, especially in random and real-world networks. There has been an immense development in the random graph literature in the recent years. We refer to the recent monographs by \cite{vdH16, van2024random, Dur25,newman2018networks} for recent results in the area.

An important of the voter model is its deep connection to random walks on graphs. In fact, the dynamics of opinion spread can be understood through a dual process: a system of \emph{coalescing random walks} evolving backwards in time (see \cite{Lig85} and also discussed later in Section~\ref{sec:VoterDef}).  As a consequence, the consensus time $\tau^u_{\rm cons}$ is closely tied to the time it takes for these independent random walks, one from each vertex, to fully coalesce.

This duality provides a powerful tool for analyzing the voter model, particularly on graphs where random walk behavior is tractable. For example, on the complete graph, the voter model reduces to the Moran model which has a simple dual structure, allowing for an explicit computation of consensus times and opinion fluctuations. Similarly, in finite-dimensional discrete tori, the coalescent process corresponds to random walks in Euclidean space, and the consensus time can be characterized via scaling limits \cite{Cox89}. On random $d$-regular graphs, the locally tree-like structure enables a detailed analysis of coalescence times using Green’s function methods, leading to precise asymptotic results \cite{CF05}. Furthermore, a precise quasi-stationary behavior of the discordant edges giving further insights on the voter evolution before the consensus time-scale can be shown, see \cite{ABHHQ22}. Recently, in \cite{avena2025}, these results have been even extended to the dynamic setting where the graphs themselves evolve over time according to an edge-rewiring mechanism preserving the regular degree structure.

However, extending these methods to more heterogeneous geometries, such as graphs with heavy-tailed degree distributions or directed structure poses significant challenges. In such settings, coalescence behavior becomes more complex due to non-uniform stationary distributions, long-range dependencies, and lack of symmetry. Understanding consensus dynamics in these networks requires new tools and various approximations. One of the major aim of the present work is to predict the behaviour of consensus time in heterogeneous directed networks with heavy-tailed degree distributions. 
\paragraph{\bf Order of mean consensus time on \emph{undirected} heterogeneous graphs.} 

Much attention has also focused on \emph{heterogeneous undirected networks} with heavy-tailed degree distributions, motivated by real-world social and information networks. In this context, statistical physics work predicts dramatically different scaling for consensus times. \cite{SAR08} studied the voter model on random graphs with an arbitrary (uncorrelated) degree distribution, showing that the mean consensus time satisfies
\begin{equation} \label{eq:Sood prediction consensus}
    \mathbf{E}[\tau^u_{\rm cons}] = \Theta\Big(H(u)\,n\,\frac{m^2_1}{m_2}\Big)\,,
\end{equation}
where $\Theta$ stands for the standard asymptotic equivalence in order, see Section \ref{Notation}, 
\begin{equation}
    m_k := \frac{1}{n} \sum_{x\in V } d^k_x, \quad k=1,2,
\end{equation}
denotes the $k$-th moment of the degree sequence,
and 
\begin{equation} \label{eq: entropy function H}
    H(u) = -(1-u)\log(1-u) - u\log(u)
\end{equation}
is an entropy factor related to the fact that the opinions are initialized as independent Bernoulli's random variables of fixed density $u\in (0, 1)$.
In particular, if the degree distribution has a power-law tail $P(k)\propto k^{-\alpha-1}$, then the average consensus time scales as $n$ for $\alpha>2$, as $n/\ln n$ for $\alpha=2$, as $n^{2-\alpha/\alpha}$ for $1<\alpha<2$, as $(\ln n)^2$ for $\alpha=1$, and remains $O(1)$ for $\alpha<1$. These predictions agree with simulations on uncorrelated scale-free networks and illustrate how extreme heterogeneity can dramatically shorten consensus times. However, mathematical proofs in these heavy-tailed regimes have been lacking. Only recently have rigorous results begun to appear. For example, \cite{FO22} studied the voter model on \emph{subcritical scale-free} inhomogeneous random graphs with independent edge connections  and showed that the time to consensus exhibits a rich phase diagram governed by different parts of the graph structure. In the directed setting, work is even less. \cite{CCPQ21} analyzed the directed configuration model with heavy-tailed in-degrees, establishing some properties of the stationary distribution and random-walk mixing time. But consensus time in directed heterogeneous graphs remain poorly understood. When the degrees are bounded, the consensus time was studied in \cite{ACHQ23} and this will form a basis of our work here also. We also mention that dynamic random graph models, in which the voter evolution influences the evolution of the underlying graph, have recently attracted attention. Various rewiring rules coupled with voter dynamics were proposed in \cite{holme2006nonequilibrium, durrett2012graph}, with predictions based on numerical simulations. Recent works such as \cite{basu2017evolving, athreya2025co, baldassarri2024} have rigorously studied these models in the setting of dense random graphs, while \cite{avena2025} considers the sparse regime. We do not consider such dynamic voter models in this article.

\subsection{ Our main contribution: Exact first-order  asymptotics of consensus time in \em{ directed heterogeneous} networks.} 

Inspired by the detailed results for the $d$‑regular ensemble and by recent advances in the analysis of directed locally tree‑like random graphs (see \cite{BCS19, CCPQ21}), the work of \cite{ACHQ23} is the first, to our knowledge, to obtain a precise asymptotic for the voter‐model consensus time in a non‑regular directed setting.  The study in \cite{ACHQ23} show that the order of the consensus time is $n$ and the precise explicit pre-constant $\vartheta_n(\mathbf{d}^+, \mathbf{d}^-)$,is expressible in terms of the in-degree and out-degree sequences, respectively, $\mathbf{d}^-$ and $\mathbf{d}^+$. 

Building on this and on the prediction in \cite{SAR08} for undirected graphs, in this paper we extend the analysis to general directed graphs whose degree sequences are sampled from a Pareto distribution with parameter $\alpha>0$ (see \eqref{eq: Pareto}).  Our main result shows that for all $\alpha>0$, including those beyond the bounded‐degree setting, the same first‐order formula holds. That is,

\begin{equation}
  \label{DigraphsForm}
  \mathbf{E}_u[\tau^u_{\rm cons}]
  \;\sim\;
  \vartheta_n(\mathbf{d}^+, \mathbf{d}^-)\,H(u)\,n.
\end{equation}
\noindent
This prediction is consistent with the undirected case in \eqref{eq:Sood prediction consensus}. In fact, as \(n\to\infty\),
\[
\vartheta_n(\mathbf{d}^+, \mathbf{d}^-)
= \Theta\bigl((m_1^-)^2/m_2^-\bigr),
\]
but our formula is more precise since it includes the explicit prefactor. This precision allows us to quantify exactly how degree heterogeneity accelerates or delays consensus. We emphasize that both the mathematical form and the resulting constants differ substantially between directed and undirected models (see \cite{ACHQ23} for a comparison). In particular, determining the analogous pre–constant in the Sood–Redner formula \eqref{eq:Sood prediction consensus} for the undirected voter model in the heterogeneous case remains an open problem.

\subsection{ Other results: Mean-field theory in directed and undirected heavy-tailed ensembles.} 
The sparse random graph models with bounded degrees studied in \cite{CF05,CCC16, ACHQ23, ABHHQ22, Cap25} satisfy mean-field conditions, which guarantee mixing and expansion properties analogous to those of the complete graph. These conditions were formalized in \cite{Oli12, Oli13, CCC16} and can be briefly summarized as follows. If the underlying graph sequence has a random-walk invariant distribution that is not too concentrated and a suitably small mixing time, as in the complete graph, then the distribution of the consensus time can be related to the expected meeting time of two independent random walks started from stationarity. In particular, \cite{Oli12,Oli13} shows that, in both directed and undirected settings, the consensus time, when rescaled by this meeting time, converges in distribution to an infinite sum of independent exponential random variables. This same limit also describes the full coalescence time of the dual random-walk system and is captured by Kingman’s coalescent. Moreover, under these mean-field conditions, \cite{CCC16} proves that the density of one opinion type evolves diffusively on the meeting-time scale and converges to the Wright–Fisher diffusion.

For the directed Pareto random-graph ensemble used to test prediction \eqref{DigraphsForm}, and for the corresponding undirected ensemble studied in \cite{SAR08}, current rigorous results confirm the mean-field conditions only in those heavy-tail regimes that guarantee sufficiently high finite moments of the degree sequence. In Section \ref{sec: MFtheory: Coalescence}, we numerically investigate which Pareto tail parameter ($\alpha>0$) satisfy these conditions in both the directed and undirected settings. We find that mean-field behavior emerges as soon as the Pareto distribution has a finite mean degree, extending beyond the regimes covered by existing proofs. While the asymptotic formulas \eqref{eq:Sood prediction consensus} and \eqref{DigraphsForm} appear universal across all heavy-tail regimes, the mean-field approximations themselves are not. For example, on the star graph with $n$ leaves, the Sood–Redner prediction still holds, yet neither the mean-field conditions nor the Kingman–coalescent and Wright–Fisher diffusion approximations can be applied. 

\subsection{Asymptotic notation}\label{Notation} We adopt the usual Bachmann–Landau notation. Given two functions \( f(n) \) and \( g(n) \): $\mathbb{N}\to\mathbb{R}^+$, we define: \( f(n) = O(g(n)) \) if there exist constants \( c, n_0 > 0 \) such that \( |f(n)| \leq c |g(n)| \) for all \( n \geq n_0 \), \( f(n) = \Theta(g(n)) \) if \( f(n) = O(g(n)) \) and \( g(n) = O(f(n)) \) ,  \( f(n) = o(g(n)) \) if \( f(n) / g(n) \to 0 \) as \( n \to \infty \),  \( f(n) = \omega(g(n)) \) if \( f(n) / g(n) \to \infty \) as \( n \to \infty \).

\subsection{Paper organization}\label{Org} 
In Section~\ref{Models}, we define the directed and undirected configuration models and describe the voter model dynamics. In Section~\ref{DirectedConsesus}, we derive asymptotic formulas for the consensus time in both settings and support our results with simulations. Section~\ref{sec: MFtheory: Coalescence} explores the mean-field approximation and its connection to coalescing random walks, while Section~\ref{sec: MFtheory: diffusions} investigates when the Wright-Fisher diffusion accurately captures the opinion dynamics. 


\section{Models}\label{Models}

\subsection{Voter model and random walk}\label{sec:VoterDef}
Consider $G=(V,E)$ to be a possibly directed, (strongly) connected graph, with vertex set $V$ and edge set $E$. We define the voter model on $G$ as the continuous-time Markov process $(\eta_t)_{t\geq 0}$ with state space $\{0,1\}^V$ and infinitesimal generator $\mathcal{L}_{vm}$ as
\begin{equation} \label{eq:generator VM}
    (\mathcal{L}_{vm}f)(\eta) = \sum_{x\in V} \sum_{y\in V} q(x,y) \,\big[f(\eta^{x\to y}) - f(\eta)\big]\,,
\end{equation}
where $ f:\{0,1\}^V \to \mathbb{R}\,$ and
\begin{equation} \label{eq:VM rates}
    q(x,y) = \frac{A(x,y)}{d^+_x}\,,
\end{equation}
while $A(\cdot,\cdot)$ is the adjacency matrix of $G$ defined as $A(x,y)\coloneqq|\{e\in E\mid e=(x,y)\}|=\#\text{ edges from }x\text{ to } y$, and $d^+_x$ denotes the number of vertices connected to $x$ via an outgoing edge. In case of an undirected graph we consider $d_x$, the usual degree of vertex $x$ in place of $d^+_x$. Finally
\begin{equation*}
    \eta^{x\to y} (z) = 
    \begin{cases}
        \eta(y),& \quad \text{if } z=x\,,\\
        \eta(x),& \quad \text{otherwise}\,. \\
    \end{cases}
\end{equation*}
For any $x\in V$ and $t\in \mathbb{R}^+$, $\eta_t(x)$ represents the state of node $x$ at time $t$ in terms of the binary state $\{0,1\}$, to be interpreted as the opinion of the individual $x$ at time $t$. In other words, the process captures the evolution of the opinion dynamics starting from the initial configuration of opinions given by $\eta_0 = \{\eta_0(x) \mid x\in V\}$. The Markovian evolution defined by the generator in \eqref{eq:generator VM} can be described as follows: Give to each vertex $x$ an exponential clock of rate $1$. When the clock associated to a vertex $x$ rings, a neighbor $y$ of $x$ is chosen uniformly at random, then vertex $x$ adopts the opinion of vertex $y$. Similarly to other interacting particle systems, such description gives rise to the so-called \emph{graphical representation} for the voter model. We refer to \cite{Lig85} and \cite{Lig99} for the details regarding graphical representation. We denote by $\mathbf{P}_u$, $u\in[0,1]$, the law of the voter model $(\eta_t)_{t\geq 0}$ with initial distribution $\eta_0 = \text{Bern}(u)^{\otimes V}$, and $\mathbf{E}_u$ its expectation.

Notice that such Markov process has two absorbing states, corresponding to the monochromatic configurations $\bar{0}$ and $\bar{1}$ consisting of all $0$'s and $1$'s, respectively. If we assume $G$ to be finite, then it can be shown that almost surely the process will reach one of the two absorbing states in finite time. This setting naturally leads to the question of determining the time such that the system reaches the absorbing states, called \emph{consensus time}, and formally defined as
\begin{equation} \label{eq: consenus time def}
    \tau^u_{\rm cons} = \inf\{t\geq0 : \eta_t \in \{\bar{0}, \bar{1}\}\}\,.
\end{equation}
One of the aims of this paper is to investigate the main properties of $\tau^u_{\rm cons}$ for specific classes of random graphs introduced in the following sections.

It is well-known that the voter model is closely related to its dual process of \emph{coalescing random walks} (\cite{Dur07, AF02}), that is, a system of independent random walks that merge into each other each time they happen to be on the same vertex of the graph. More precisely, let  $(X_t)_{t\ge0}$ denote the \emph{continuous-time random walk} on $G=(V,E)$, which is the Markov process with state space $V$  and infinitesimal generator given by
\begin{equation}\label{CTRW}
	\mathcal{L}_{\rm rw}f(x)=\sum_{y\in V}q(x,y)\,\left[f(y)-f(x)\right]\,,\qquad f:V\to\R\,,
\end{equation}
with $q(\cdot,\cdot)$ as in \eqref{eq:VM rates}. 
In the following we will assume $G$ to be strongly connected, that is, for each pair of vertices there exists a finite path connecting them. The latter requirement immediately implies that the random walk on $G$ admits a unique stationary distribution, which we denote by $\pi$.
Finally, we consider the standard \emph{mixing time} $t_{\rm mix}$, i.e. the first time such that the distribution of the random walk on the graph $G$ is close enough to $\pi$. 


\subsection{Configuration models} \label{subsection:CM defs}

The classical (undirected) Configuration Model introduced by Bollob\'as \cite{Bol80} in the early '80s and can be defined as follows. See also \cite{vdH16} for a modern introduction to the topic. For any $n\in\mathbb{N}$, let $[n]:= \{1,\dots,n\}$ be a set of $n$ labeled nodes. For any vertex $x\in[n]$, let $d_x$ be its degree, that is the number of edges that are connected to $x$. Define $\mathbf{d}_n = (d_1,\dots,d_n)$ to be a \emph{deterministic} degree sequence with the following constraint
\begin{equation} \label{graphical-undirected}
 \sum_{x\in[n]}d_x = 2\,\ell\,,
\end{equation}
for some $\ell=\ell_n$. The randomness of the model comes from the mechanism in which the edges are formed. This is a result of the following uniform pairing procedure involving \emph{stubs}, i.e. half-edges. Assign to each vertex $x\in[n]$, $d_x$ labeled stubs. At each step, select a stub $e$ that was not matched in a previous step (the order of the selection of $e$ is irrelevant), and a uniform at random stub $f$ among the unmatched ones. Then match them and add the edge $ef$ between the vertex incident to $e$ and the one incident to $f$ to the edge set $E$. Continue until there are no more unmatched stubs. Note that this is possible given the assumption \eqref{graphical-undirected} that the sum of the stubs is even. This random procedure gives rise to a so-called \emph{configuration}, and it uniquely determines the corresponding random graph $G=G_n=([n],E)$. The graph may have multiple edges and (multiple) self-loops. Importantly, this means that the number of vertices connected to a vertex is in general different from (and typically smaller than) the number of edges connected to the same vertex (the degree). We say that a graph $G_n$ is sampled from the Configuration Model CM $=$ CM($\mathbf{d}_n$) with a given degree sequence $\mathbf{d}_n$, if it is sampled according to the procedure above.

The Directed Configuration Model is a natural generalization to the directed setting of the classical version of the model defined in the previous paragraph. We give its definition for completeness. Fix $n\in\mathbb{N}$, and let $[n]$ be the vertex set as above. Similarly to the undirected counterpart, let $\mathbf{d}^+=\mathbf{d}^+_n=(d^+_x)_{x\in[n]}\in \mathbb{N}^n$ and $\mathbf{d}^-=\mathbf{d}^-_n=(d^-_x)_{x\in[n]}\in \mathbb{N}_0^n$ be two finite \emph{deterministic} sequences such that 
\begin{equation} \label{graphical}
    m=m_n\coloneqq\sum_{x\in[n]}d_x^+ = \sum_{x\in[n]}d_x^-\, .
\end{equation}
For any vertex $x\in[n]$, let $d_x^+$ (resp. $d_x^-$) be its out-degree (resp. in-degree), that is the number of directed edges that are exiting (resp. entering) $x$. At this point we perform a pairing procedure on the same line as the undirected one, but this time we need to match in-stubs, called \emph{heads}, only with out-stubs, \emph{tails}. At each step, select a tail $e$ that was not matched in a previous step, and a uniformly random head $f$ among the unmatched ones, then match them and add the directed edge $ef$ from the vertex incident to $e$ to the one incident to $f$ to the edge set $E$. Continue until there are no more unmatched heads and tails. As before, note that the constraint in \eqref{graphical} ensures that such uniform matching ends without any stub left unmatched. 
Again, note that there can be multiple edges and (multiple) self-loops, which means that the numbers of incoming and outgoing edges of a vertex do not in general coincide with the numbers of vertices that that vertex is connected to in the incoming and outgoing direction, respectively.
We say that a graph $G$ is sampled from the Directed Configuration Model DCM $=$ DCM($\mathbf{d}^+_n,\mathbf{d}^-_n$) with a given bi-degree sequence $\mathbf{d}^+_n,\mathbf{d}^-_n$, if it is sampled according to the procedure above.

Assume now that the degree sequence of the CM and DCM is sampled from a power-law distribution  
with exponent $\alpha>0$, i.e.
\begin{equation} \label{eq: Pareto}
    \mathsf{P}(D_\alpha >x) \sim x^{-\alpha}\,,\quad \text{as} \ x\to \infty\,.
\end{equation}
Without loss of generality, we can take $D_\alpha\sim \text{Pareto}(\alpha)$ and denote its density by $f_{D_\alpha}$. In particular, we will consider a Pareto distribution with exponent $\alpha$ and minimal value $x_{\rm min}\geq1$,
so its density will be of the type $f_{D_\alpha}(x) = \alpha x^{\alpha}_{\rm min}\,x^{-\alpha-1}$, $x>0$. For any $\alpha>0$, we define as $\alpha$-CM (resp. $\alpha$-DCM) the CM (resp. DCM) having degree sequence $\mathbf{d}_n$ given by i.i.d. of $D_\alpha$, and denote by $\mathsf{P}_\alpha$ and $\mathsf{E}_\alpha$ the law and the expectation of $\mathbf{d}_n$. We performed several simulations involving $\alpha$-DCM ensembles. The in- and out-degree sequences were sampled as follows: first, we sampled the in-degrees according to $D_\alpha$, and then we set the out-degrees as a random permutation of the in-degrees. This approach ensured that both distributions were correctly sampled and that the graphical assumption \eqref{graphical} was always satisfied.

The rate of the largest degree plays a crucial role in the asymptotics we describe later. Given an i.i.d. sequence $\{d_1,\dots,d_n\}$ such that $d_x \sim D_\alpha$ for all $x\in[n]$, it holds that 
\begin{equation} \label{eq: d_max asymptotic}
    \lim_{n\to \infty}\mathsf{P}(d_{\rm max}/n^{1/\alpha}\le y) = \mathsf{P}(Y\le y), \quad y>0\,,
\end{equation}
where $Y\sim \text{Fr\'echet}(\alpha)$ and hence note that $d_{\max}= \Theta_\mathsf{P}(n^{1/\alpha})$.
We use a slight abuse of notation by denoting by $\P$ the law of $G$, both for the CM and the DCM, and we will specify it whenever it will not be clear from the context. We will be interested in studying the asymptotic regime in which $n\to\infty$, and we will say that $G$ has a certain property $E_n$ \emph{with high probability (w.h.p.)}, if
\begin{equation} \label{with high probability}
\P(E_n) \to 1\,,  \text{ as } n\to\infty\,.    
\end{equation}

\subsubsection{\bf Stationary measure and mixing}  In later sections, we will be interested in the known results for the stationary distribution $\pi$ and the mixing time $t_{\rm mix}$ defined in Section \ref{sec:VoterDef}, for the random walk on the CM and DCM. Note that, as for any undirected graph, the stationary distribution $\pi$ of the random walk on CM is given by the deterministic measure proportional to the degrees,
despite the graph being random. Whereas for the DCM ensemble, the measure $\pi$ is an unknown \emph{random measure} that may strongly depend on the realization of the random graph (with the exception of the Eulerian case in which $d_x^+=d_x^-$ for all $x$, where the expression of $\pi$ goes back to the undirected case). A complete description of the distribution of $\pi$ is unknown. Nevertheless, in recent years some crucial properties were derived under different assumptions on the degree sequence (\cite{CCPQ21,CP20,CQ20,CQ21b}). 

The mixing time is a random variable with respect to the graph ensemble both in the CM and in the DCM. In the CM setting, assuming that the minimal degree is at least $3$, it is expected that \cite[Theorem 6.4.1]{Dur25}, 
\begin{equation}
    t_{\rm mix} = \mathcal{O}(\log(n)) \quad w.h.p.
\end{equation}
For the DCM the known results are very promising, but less complete. Indeed, more precise results, the so called cutoff at entropic time, have been proven for some regimes of the bi-degree sequence in \cite{BCS18,BCS19,CCPQ21}.




\section{Consensus time for CM and DCM}\label{DirectedConsesus}
 We begin this section by recalling the prediction in \cite{SAR08} for the expected consensus time for independent-edges (uncorrelated), weighted random graphs, having transition matrix replaced by its expectation, that is
\begin{equation}
    A(x,y) = A(y,x) = \frac{w_x\,w_y}{\sum_{x\in[n]} w_x}\,,
\end{equation}
where $\{w_x\}_{x\in [n]}$ are weights that on average equal to the degrees $\{d_x\}_{x\in [n]}$. Consider the case when $d_x$ is sampled from the Pareto distribution $D_\alpha$ for all $x$. 

Observe first that, by the strong law of large numbers, for $n$ large the $i$-th empirical moment $\Hat{m}_i = \frac{1}{n}\sum_{x\in[n]}d^i_x$ of the degree sequence can be approximated by their truncated moment
\begin{equation}
    m_i = \int x^i \, f_{D_\alpha}(x)\,dx = \int^{d_{\rm max}}_{x_{\rm min}} x^i \, f_{D_\alpha}(x)\,dx\,.
\end{equation}
Since we are interested only in the first and second moments, let us explicitly compute the latter in the $i=\{1,2\}$ cases, exploiting the $d_{\rm max}$ asymptotic in \eqref{eq: d_max asymptotic}. In the following $C$ is a positive random variable depending only on $Y$, $x_{\rm min}$ and $\alpha$.

\[
m_1 = \int_{x_{\rm min}}^{d_{\max}} x^{-\alpha} \,dx \sim  \begin{cases}  
n^{\frac{1-\alpha}{\alpha}} - C = \mathcal{O}(1), & \alpha > 1 \\  
\log(n), & \alpha = 1 \\  
n^{\frac{1-\alpha}{\alpha}}, & 0 < \alpha < 1  
\end{cases}
\]

\[
m_2 = \int_{x_{\rm min}}^{d_{\max}} x^{1-\alpha} \,dx \sim  \begin{cases}  
n^{\frac{2-\alpha}{\alpha}} - C = \mathcal{O}(1), & \alpha > 2 \\  
\log(n), & \alpha = 2 \\  
n^{\frac{2-\alpha}{\alpha}}, & 1 < \alpha < 2 \\  
n, & \alpha = 1 \\  
n^{\frac{2-\alpha}{\alpha}}, & 0 < \alpha < 1  
\end{cases}
\]

Therefore, following the author's predictions for their mean-field model, the expected consensus time has, up to a pre-constant, the order defined in Equation \eqref{eq:Sood prediction consensus}.
We can interpret \eqref{eq:Sood prediction consensus} as a product between a \emph{volume} factor arising from the size of the largest connected component, $n$ in this case as we are assuming the ensembles to be in their connected regime, and a \emph{volatility} factor indicating how the consensus can be sped up or slowed down by arranging the degree sequence of the network. We will perform a more detailed analysis of such result in the following Section \ref{sec:refined results undirected}.

\subsection{Refined results for undirected case} \label{sec:refined results undirected}
As mentioned above the behaviour of the expected consensus time can be derived from the expressions for $m_1, m_2$. Indeed, for any $\alpha >0$ it follows that
\begin{equation} \label{eq: expected consensus formula}
     \mathbf{E}[\tau^u_{\rm cons}] \sim c\, \log^a(n)\, n^b\,,
\end{equation}
where $c$ is a (possibly random) almost surely positive constant depending only on $\alpha,u$ and not depending on $n$, while the following are the corresponding values for the parameters $a$ and $b$.
\begin{table}[H]
    \centering
    \begin{tabular}{|c|c|c|c|c|c|}
        \hline
        & \textbf{$\alpha>2$} & \textbf{$\alpha =2$} & \textbf{$1<\alpha <2$}& \textbf{$\alpha =1$} & \textbf{$0<\alpha<1$} \\ 
        \hline
        $a$ & 0 & -1 & 0 & 2 & 0\\ 
        \hline
        $b$ & 1 & 1 & $2(\alpha-1)/\alpha$ & 0 & 0\\ 
        \hline
    \end{tabular}
    \label{tab:a and b values}
\end{table}
 In particular, the infinite mean case ($\alpha <1$) leads to an expected consensus reached in constant time as $n$ increases. 
\begin{figure*}[t] 
    \centering
    \begin{subfigure}{0.45\textwidth}
        \centering
        \includegraphics[width=\linewidth]{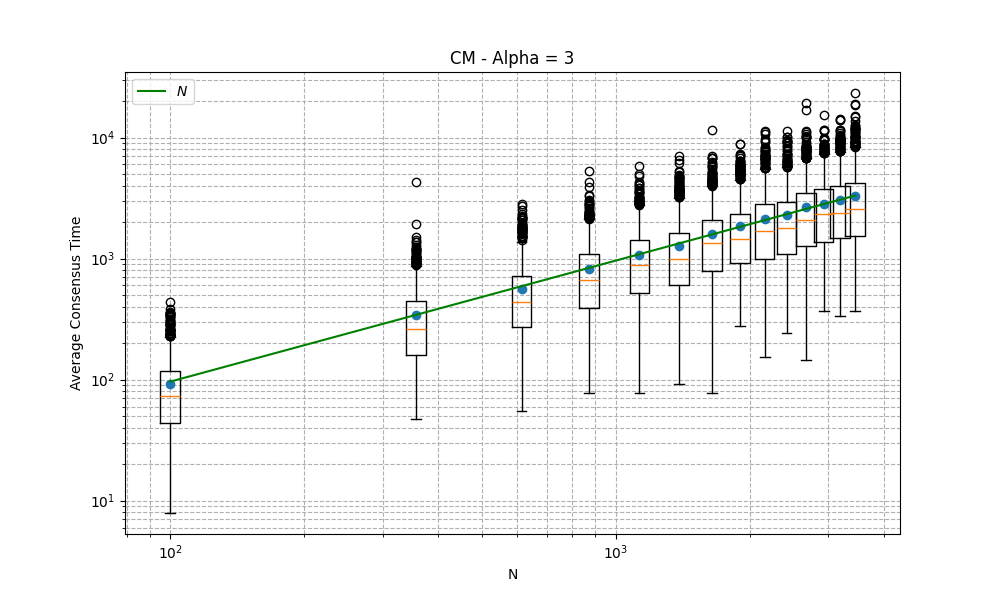}
        \caption{$\alpha = 3$, log-log scale.}

    \end{subfigure}
    \hfill
    \begin{subfigure}{0.45\textwidth}
        \centering
        \includegraphics[width=\linewidth]{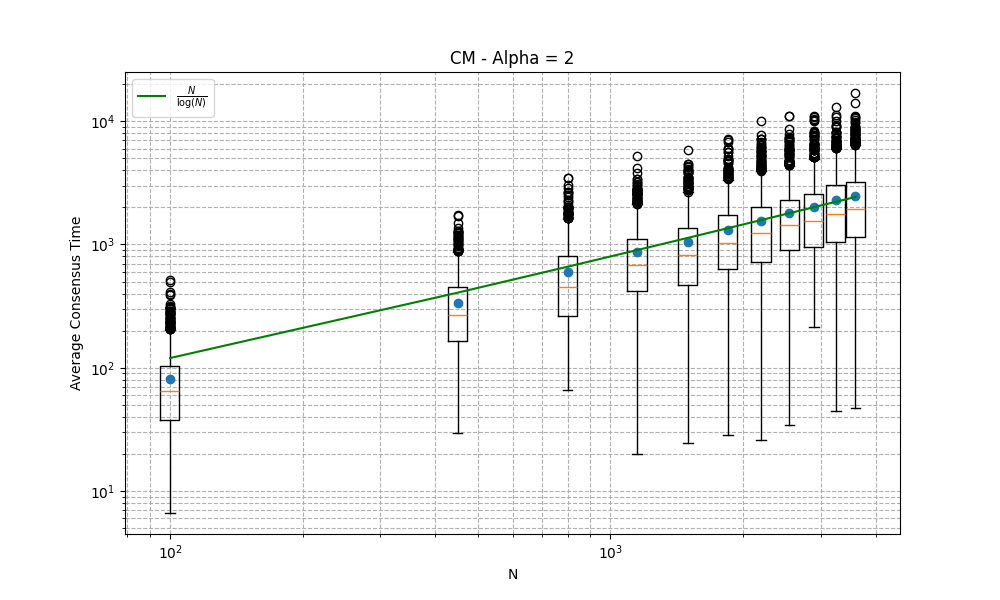}
        \caption{$\alpha = 2$, log-log scale.}

    \end{subfigure}
    
    
    \begin{subfigure}{0.45\textwidth}
    \centering
    \includegraphics[width=\linewidth]{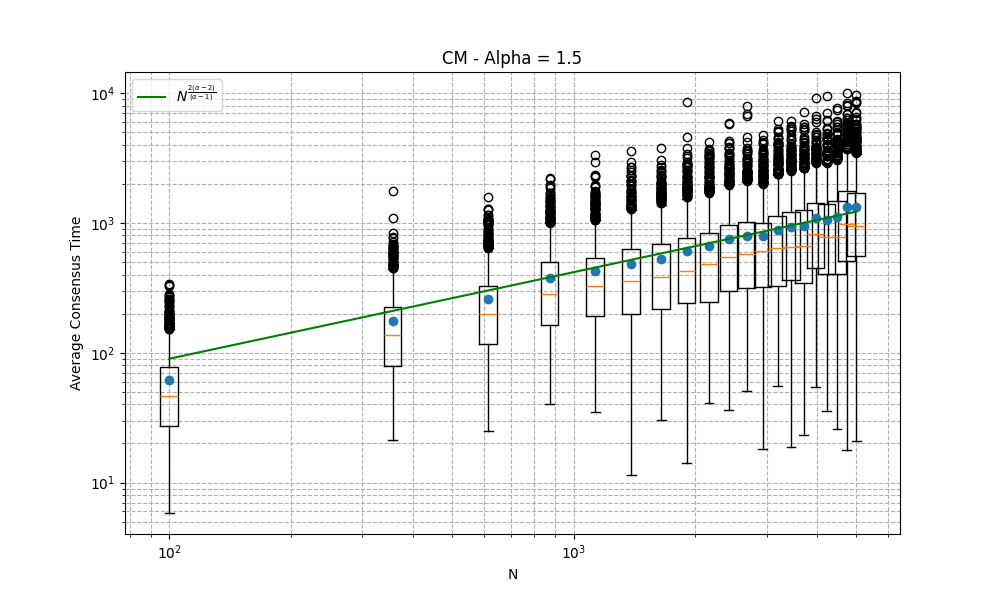}
    \caption{$\alpha = 1.5$, log-log scale.}

    \end{subfigure}
    \hfill
    \begin{subfigure}{0.45\textwidth}
        \centering
        \includegraphics[width=\linewidth]{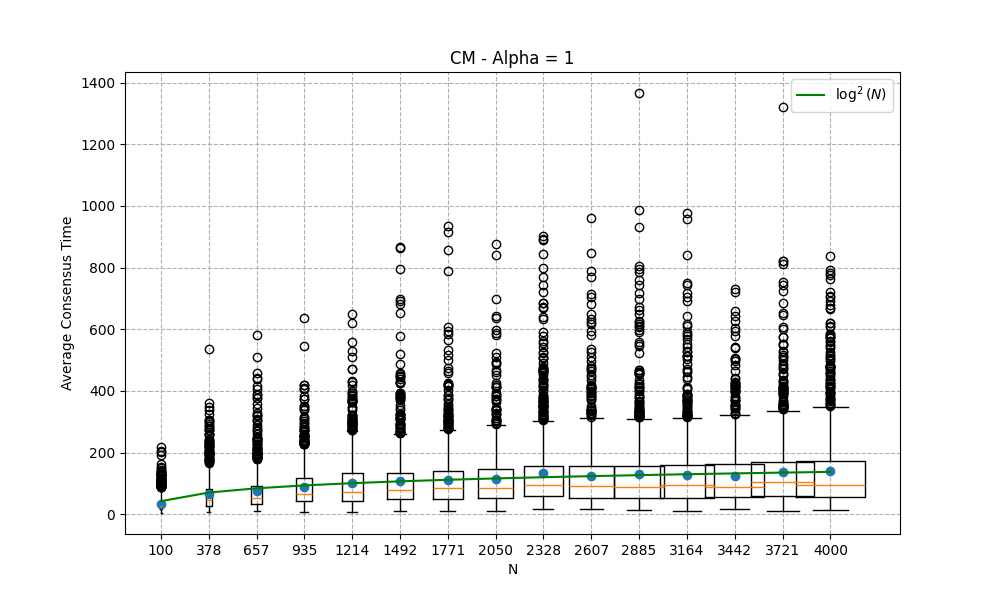}
        \caption{$\alpha = 1$}

    \end{subfigure}
    

    \begin{subfigure}{0.4\textwidth}
        \centering
        \includegraphics[width=\linewidth]{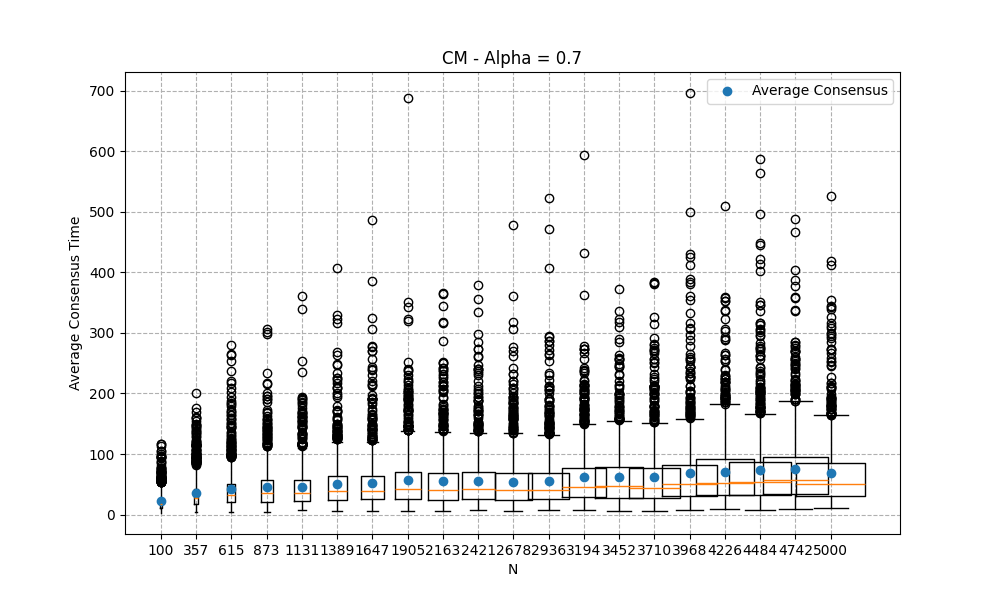}
        \caption{$\alpha = 0.7$}

    \end{subfigure}

    \caption{All figures represent the comparison between the box plots of the consensus time in the $\alpha$-CM, for different ranges of $\alpha$ and increasing values of $n$. The blue dots denotes, for each $n$, the average over the voter model, the $\alpha$-CM and the degree sequence realizations, while the box plots are taken over all the possible voter model outcomes. In figures (a)-(c) the plot is given in log-log scale and the green lines represent the predicted expected consensus time in terms of powers of $n$ as described in \eqref{eq: expected consensus formula}. In figures (d)-(e) the scale is not logarithmic as it pictures values that are slowly or not growing as $n$ increases. Here, for each $n$ ranging between 100 to 5000, we sampled 50 degree sequences, for each of which 20 graph realizations, for each of which 10 voter realizations.}
    \label{fig:CM}
\end{figure*}
In this section, we want to perform a more refined analysis regarding the study of $\mathbf{E}[\tau^u_{\rm cons}]$ for the $\alpha$-CM model defined in Subsection \ref{subsection:CM defs}.
We present our study from multiple perspectives: we confirm through simulations that the asymptotic behavior of the expected consensus time aligns with the predictions in \cite{SAR08}; we provide a better understanding of the interplay between the maximal degree and the consensus time dynamics; we analyze the results under different levels of quenching, assessing their stability with respect to various sources of randomness; we conjecture a more general result that applies across different connectivity regimes of the $\alpha$-CM.

\begin{figure}[ht]
    \centering
    \includegraphics[width=0.8\linewidth]{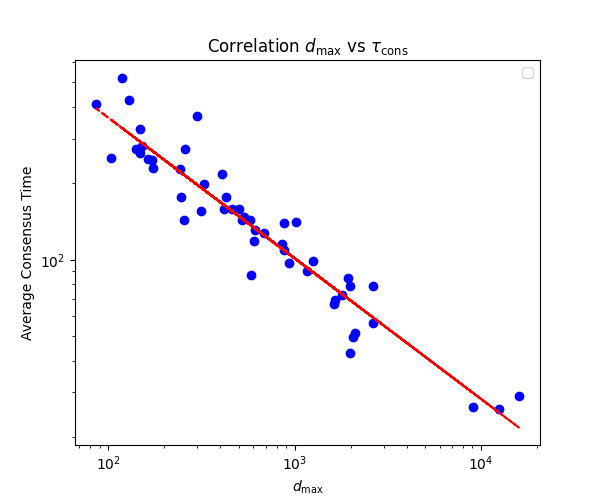}
    \caption{Scatter plot of $d_{\rm max}$ and $\mathbf{E}[\tau^u_{\rm cons}]$. Here we simulated for a fixed value of $n=1000$ and a single degree sequence, 50 graph realizations, for each of which 10 voter realizations. Each blue dot represents the average of the 10 voter runs for each of the 50 graph realizations, while the red line is just a linear regression of the outcomes.}
    \label{fig:correlation d_max-consensus}
\end{figure}

\begin{figure*}[t]
    \centering
    \begin{subfigure}{0.45\textwidth}
        \centering
        \includegraphics[width=\linewidth]{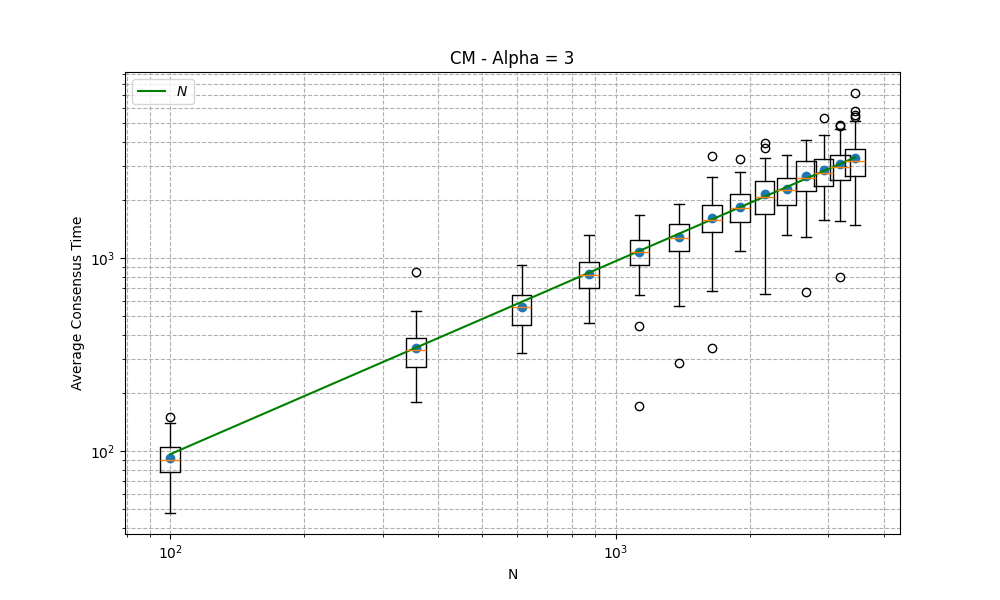}
        \caption{$\alpha = 3$, log-log scale.}

    \end{subfigure}
    \hfill
    \begin{subfigure}{0.45\textwidth}
        \centering
        \includegraphics[width=\linewidth]{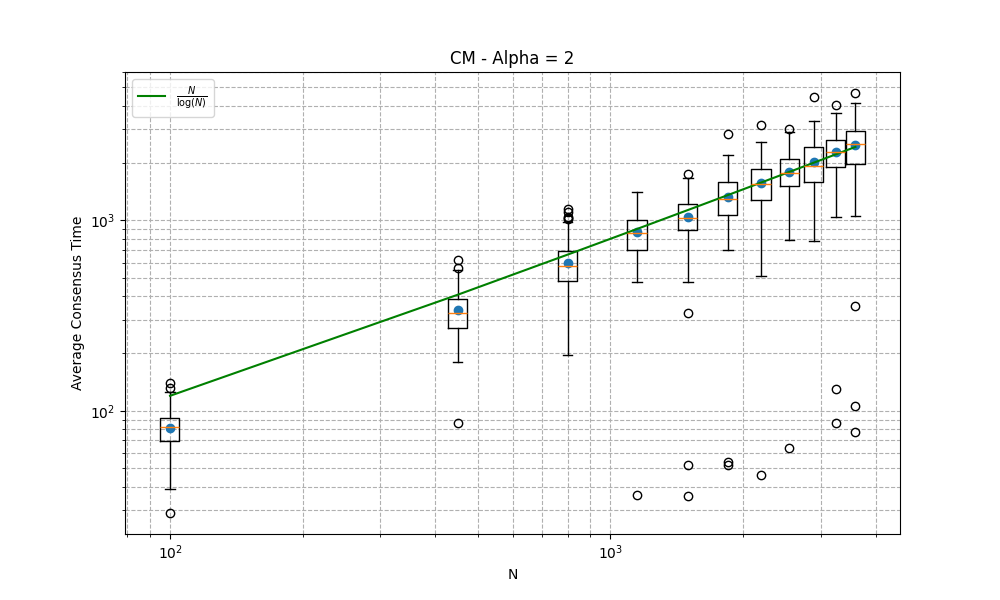}
        \caption{$\alpha = 2$, log-log scale.}

    \end{subfigure}
    
    
    \begin{subfigure}{0.45\textwidth}
    \centering
    \includegraphics[width=\linewidth]{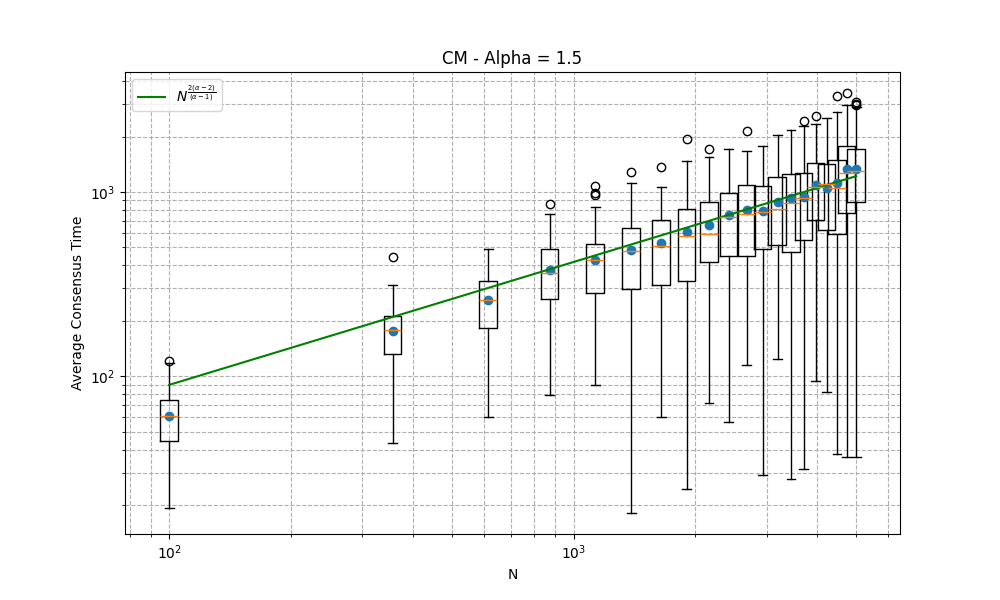}
    \caption{$\alpha = 1.5$, log-log scale.}

    \end{subfigure}
    \hfill
    \begin{subfigure}{0.45\textwidth}
        \centering
        \includegraphics[width=\linewidth]{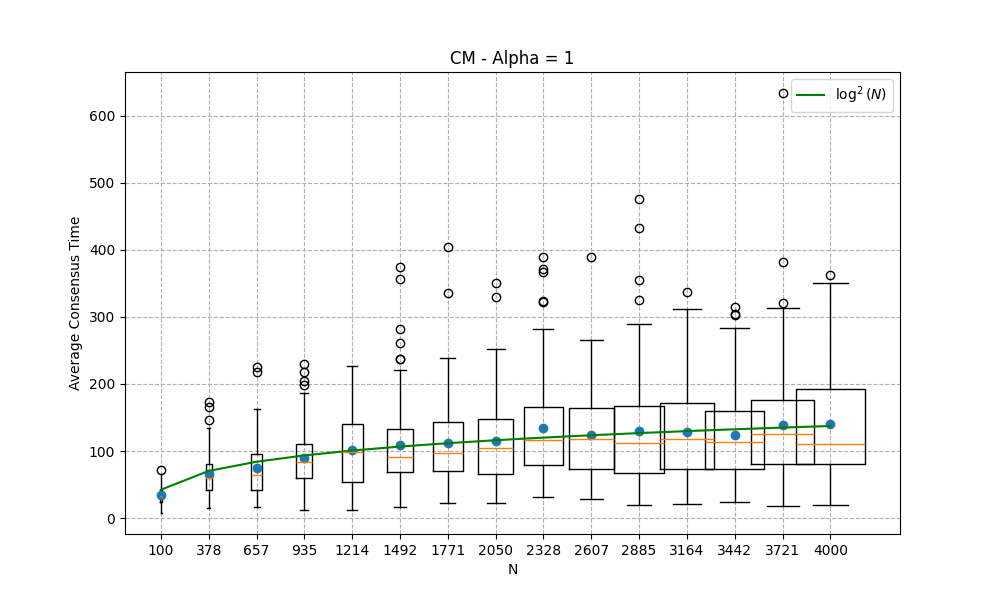}
        \caption{$\alpha = 1$}

    \end{subfigure}
    

    \begin{subfigure}{0.45\textwidth}
        \centering
        \includegraphics[width=\linewidth]{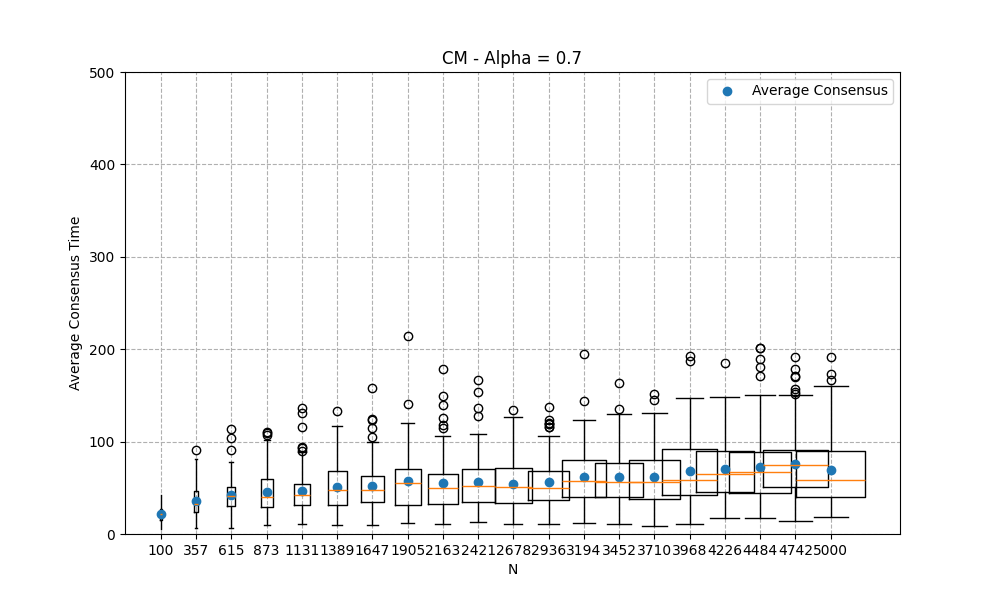}
        \caption{$\alpha = 0.7$}

    \end{subfigure}

    \caption{The same realizations as in Figure \ref{fig:CM}, where the outcomes are averaged with respect to the voter model. As a consequence, the upper skewness of the outliers is significantly reduced. Here, for each $n$ ranging between 100 to 5000, we sampled 50 degree sequences, for each of which 20 graph realizations, for each of which 10 voter realizations.} 
    \label{fig:CM_average_VM}
   
\end{figure*}

Figure \ref{fig:CM} captures the expected behavior of the $\tau^u_{\rm cons}$ over different values of the parameter $\alpha$. There are a few remarks that need to be addressed regarding the outcomes of the consensus times. The box plot visualization shows the presence of outliers, all within the upper region above the third quartile. This evidence, strengthened by the fact that the mean value (blue dot) is always above the median (yellow line), can be justified by an intrinsic skewness of the consensus time distribution, and a strong negative correlation between the consensus time and the maximal degree of the related underlying network. The first can be easily observed in Figure \ref{fig:CM_average_VM}, where we replace all consensus times outcomes by their average, for each degree sequence and graph realizations. Regarding the relationship between $d_{\rm max}$ and $\mathbf{E}[\tau^u_{\rm cons}]$, from Figure \ref{fig:correlation d_max-consensus} we deduce that they are negatively correlated. We conclude our argument by noticing that, as stated in \eqref{eq: d_max asymptotic}, $d_{\rm max}$ properly rescaled converges in probability to a $\text{Fr\'echet}(\alpha)$, which is a heavy-tailed skewed to the left random variable.

\subsubsection{Annealed vs Quenched}
 We begin by noting that $\tau^u_{\rm cons}$ involves four distinct, dependent sources of randomness.: the initial condition of the opinions, the degree sequence, the random graph and the voter model evolution.  For the $\alpha$-CM, so far we analyzed the \emph{annealed} dynamic of the process, averaging on all the layers, while we are also interested in the study of the process on different layers of \emph{quenching}. In fact, a more precise statement of \eqref{eq: expected consensus formula} is 
\begin{equation} \label{eq: annealed expected consensus formula}
         \Bar{\mathbb{E}}_\alpha\big[\mathbf{E}[\tau^u_{\rm cons}]\big] \sim c\, \log^a(n)\, n^b\,,
\end{equation}
where $\Bar{\P}_\alpha$ is the joint measure with respect the randomnesses of the degree sequence $\mathsf{P}_\alpha$ and the random graph $\P$, and $\Bar{\mathbb{E}}_\alpha$ the corresponding expectation.

\begin{figure}[ht]
    \centering
    \begin{subfigure}{0.45\textwidth}
        \centering
        \includegraphics[width=\linewidth]{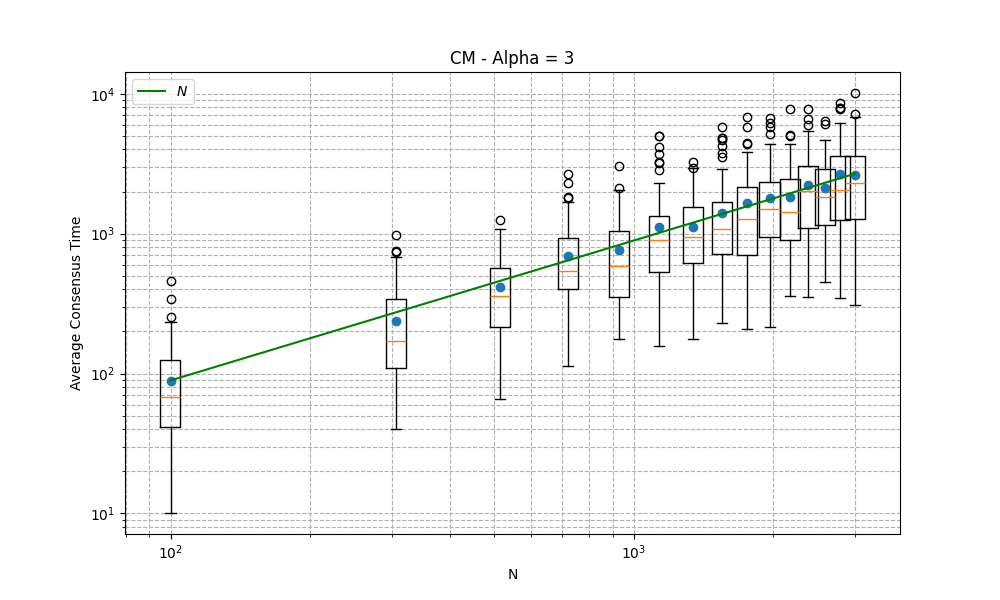}
        \caption{$\alpha = 3$, log-log scale.}
  
    \end{subfigure}
    \hfill
    \begin{subfigure}{0.45\textwidth}
        \centering
        \includegraphics[width=\linewidth]{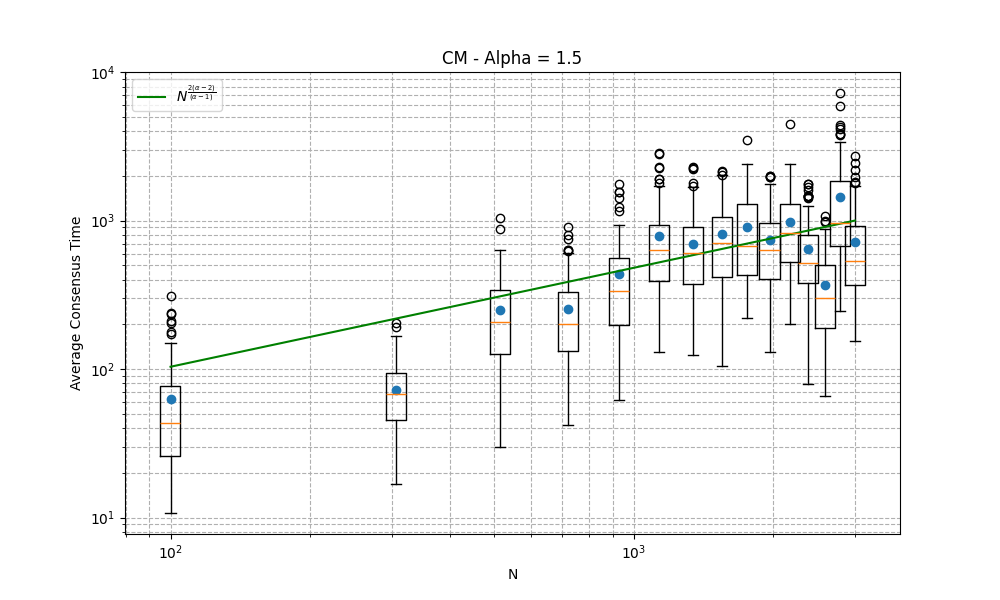}
        \caption{$\alpha = 1.5$, log-log scale.}

    \end{subfigure}

    \caption{For each $n$ ranging between 100 to 5000, we sampled one degree sequence, for which we sampled 20 graph realizations, for each of which 50 voter realizations, i.e., we quenched only on the degree sequence.} 
    \label{fig:Graph_Quenched_CM_average_VM}
\end{figure}

We implemented several experiments with the aim of exploring the stability of the random variable $\mathbf{E}[\tau^u_{\rm cons}]$ with respect to the different sources of randomness that interact within it. More precisely, we ask ourself whether it is reasonable to conjecture that a quenched result may hold, i.e. if \eqref{eq: expected consensus formula} holds with high probability with respect to $\Bar{\P}_\alpha$.
We treat only the cases $\alpha>1$, since the infinite mean regime is too unstable and has a more subtle predicted asymptotic to match with the simulated one. Based on what we can observe in Figure \ref{fig:Graph_Quenched_CM_average_VM} and in view of the results depicted in the previous section for the directed model, we predict that we can expect a fully quenched result as in Equation \eqref{eq: expected consensus formula} only for $\alpha>2$. This conclusion is in line with the analysis made in the previous section, where the variability of the graph, partially induced by the variability of the degree sequence, plays a crucial role and cannot be addressed without a precise first-order asymptotic.

\begin{figure}[ht]
    \centering
    \begin{subfigure}{0.45\textwidth}
        \centering
        \includegraphics[width=\linewidth]{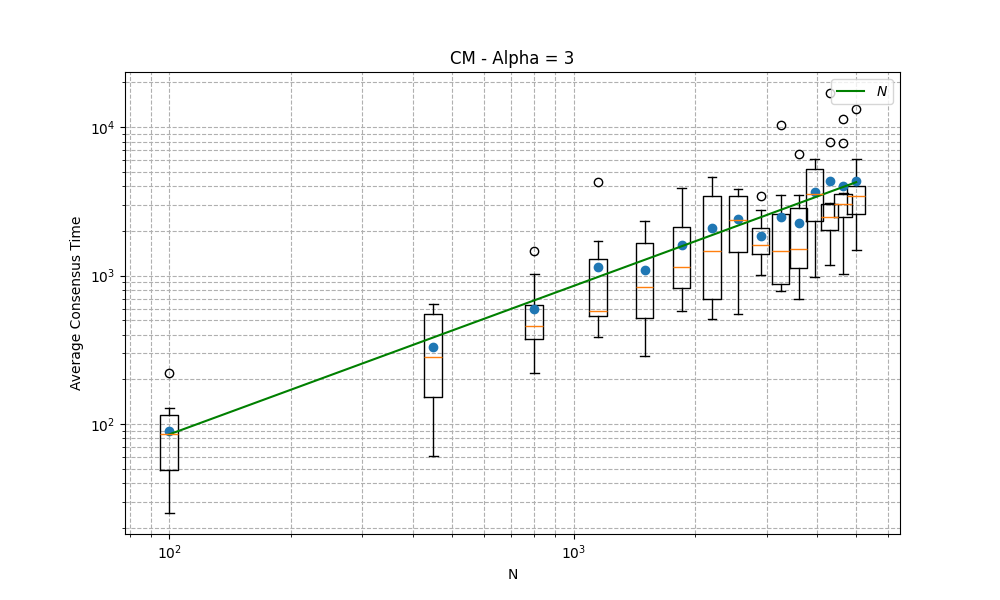}
        \caption{$\alpha = 3$, log-log scale.}
  
    \end{subfigure}
    \hfill
    \begin{subfigure}{0.45\textwidth}
        \centering
        \includegraphics[width=\linewidth]{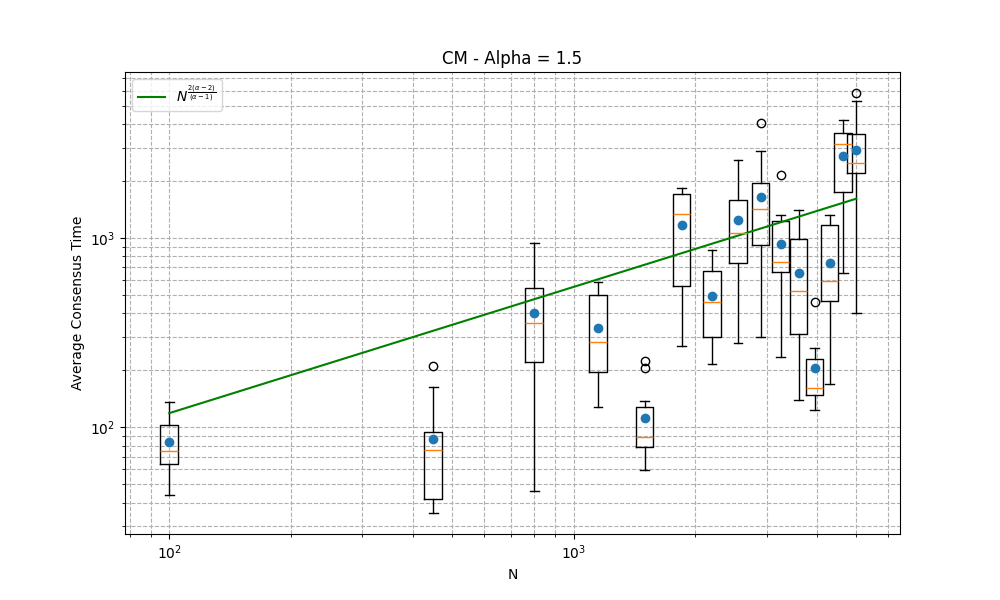}
        \caption{$\alpha = 1.5$, log-log scale.}

    \end{subfigure}

    \caption{For each $n$ ranging between 100 to 5000, we sampled one degree sequence, for which we sampled one graph realizations, for which we sampled 50 voter realizations, i.e., we quenched on both the degree sequence and the graph.} 
    \label{fig:All_Quenched_CM_average_VM}
\end{figure}

\subsection{Consensus time for directed case}
We now do a similar analysis DCM model. For a given degree sequence ($\mathbf{d}^+, \mathbf{d}^-$) define 
\begin{equation}\label{Dmoments}
    \Hat{m}_i^+ = \frac{1}{n}\sum_{x\in[n]} (d_x^+)^i \quad \text{and} \quad \Hat{m}_i^- = \frac{1}{n}\sum_{x\in[n]} (d_x^-)^i
\end{equation}
to be the $i$-th moments of the out-degree (resp. in-degree) sequence.
We now introduce some relevant functions of the degree sequence of the DCM.
Let
\begin{equation}\label{eq:def-rho-gamma-delta-beta}
\begin{split}
\delta&=\frac{m}{n}\,,\quad \qquad\qquad\beta=\frac1{m}\sum_{x\in[n]}(d_x^-)^2\,, \\ \rho&=\frac{1}{m}\sum_{x\in[n]}\frac{d_x^-}{d^+_x}
\,,\quad\gamma=\frac{1}{m}\sum_{x\in[n]}\,\frac{(d^-_x)^2}{d^+_x}\,,
\end{split}
\end{equation}
where $m=n\Hat{m}_1^-=n\Hat{m}_1^+$ as in \eqref{graphical}. Notice that $\beta$ and $\delta$ equal to $\Hat{m}_2^-/\Hat{m}_1^-$ and $\Hat{m}_1^-$, respectively.
Moreover, define
\begin{equation}\label{theta}
    \vartheta= \vartheta_n( \mathbf{d}^+, \mathbf{d}^-)\coloneqq \frac{\delta}{\frac{\gamma-\rho}{1-\rho}\,\Big(1- \frac{1-\sqrt{1-\rho}}{\rho}\Big)+\beta-1}\,.
\end{equation}
We conjecture that for \emph{any} directed graph sampled from DCM such that the above quantities are well-defined, we get with high probability
\begin{equation} \label{eq: conjecture genereal expected consensus}
   |\mathbf{E}[\tau^u_{\rm cons}] - H(u)\,\vartheta\,n | \underset{n\to\infty}{\longrightarrow} 0\,,
\end{equation}
where $H(\cdot)$ as in \eqref{eq: entropy function H}.
In particular, if we consider the $\alpha$-DCM setting defined in Subsection \ref{subsection:CM defs}, for any $\alpha>0$ we get the same asymptotic as in \eqref{eq: conjecture genereal expected consensus} where the convergence holds in probability with respect to $\Bar{\mathbb{P}}_\alpha$, the joint measure with respect the randomnesses of the degree sequence and the random graph.

The voter model on DCM was first studied in \cite{ACHQ23}, under the assumption of having a deterministic degree sequence $\mathbf{d_n}$ with \emph{bounded degrees}, i.e. $d_{\max}^\pm = \mathcal{O}(1)$.
In their main result, the authors prove the convergence of the distribution of the consensus time and the first-order asymptotic of its expectation. In particular, they show the validity of \eqref{eq: conjecture genereal expected consensus} with high probability with respect to $\P$. 

We show via simulations the validity of the conjecture \eqref{eq: conjecture genereal expected consensus} in the $\alpha$-DCM case. Informally, we claim that for any $\alpha>0$ the expected consensus time on $\alpha$-DCM has the same asymptotic above, where $\vartheta$ is no longer uniformly bounded but possibly $n$ dependent. We start our analysis with the result given in \cite{ACHQ23}, observing that under bounded degrees assumption $\vartheta = \Theta(1)$, leading to a linear expected consensus. Recall from \eqref{eq: d_max asymptotic} that, in our setting where $d_x^\pm \sim \text{Pareto}(\alpha)$, we have $d_{\rm max}^\pm = \Theta_{\mathsf{P}_\alpha}\left(n^{1/\alpha}\right)$ for any $\alpha>0$, so the uniform boundedness of the maximal in- and out-degrees no longer holds. The assumption of bounded degree is reflected in a light-tailed distribution, that is, $\alpha= \infty$.


\begin{figure*}[t] 
    \centering
        \begin{subfigure}{0.45\textwidth}
        \centering
        \includegraphics[width=\linewidth]{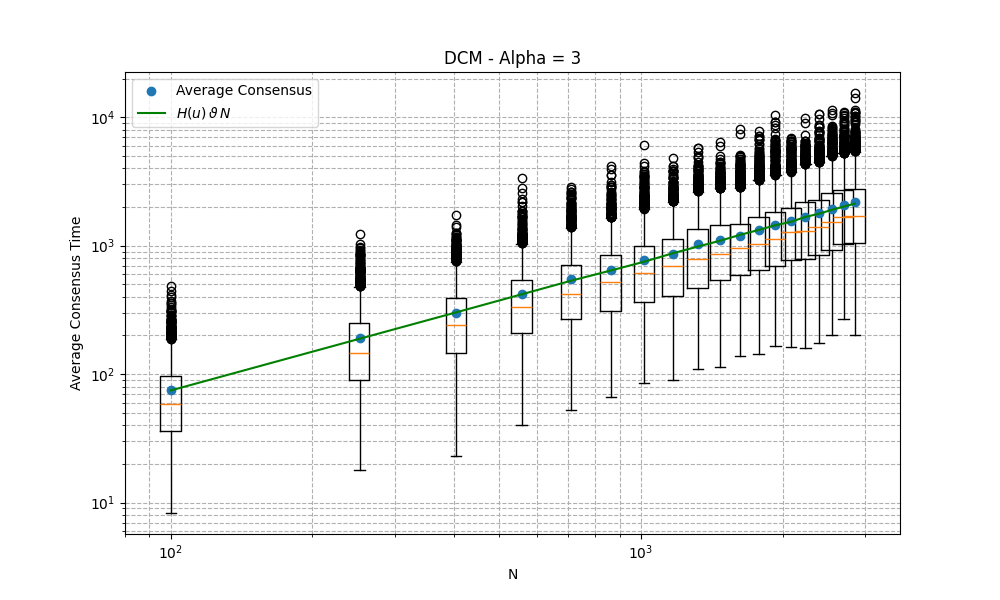}
        \caption{$\alpha = 3$, log-log scale}

    \end{subfigure}
    \hfill
    \begin{subfigure}{0.45\textwidth}
        \centering
        \includegraphics[width=\linewidth]{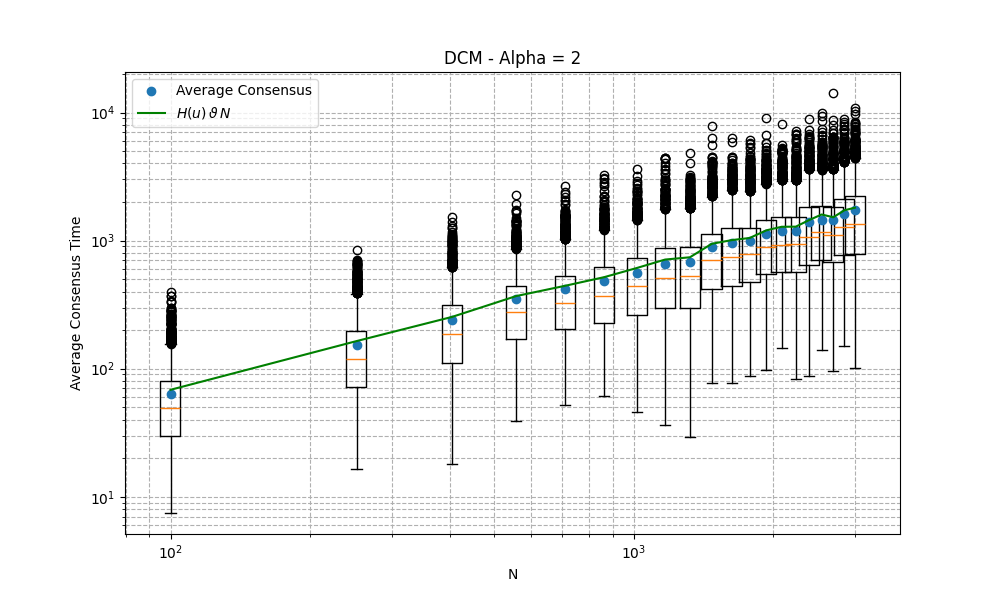}
        \caption{$\alpha = 2$, log-log scale.}

    \end{subfigure}
    
    \vspace{0.5cm}  
    \begin{subfigure}{0.45\textwidth}
        \centering
        \includegraphics[width=\linewidth]{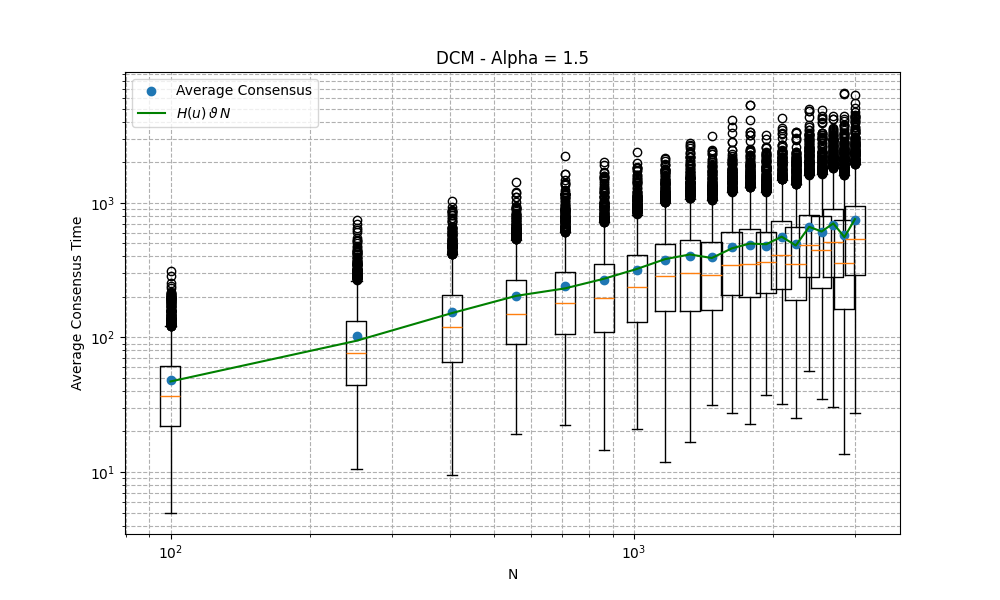}
        \caption{$\alpha = 1.5$, log-log scale.}

    \end{subfigure}
    \hfill
    \begin{subfigure}{0.45\textwidth}
        \centering
        \includegraphics[width=\linewidth]{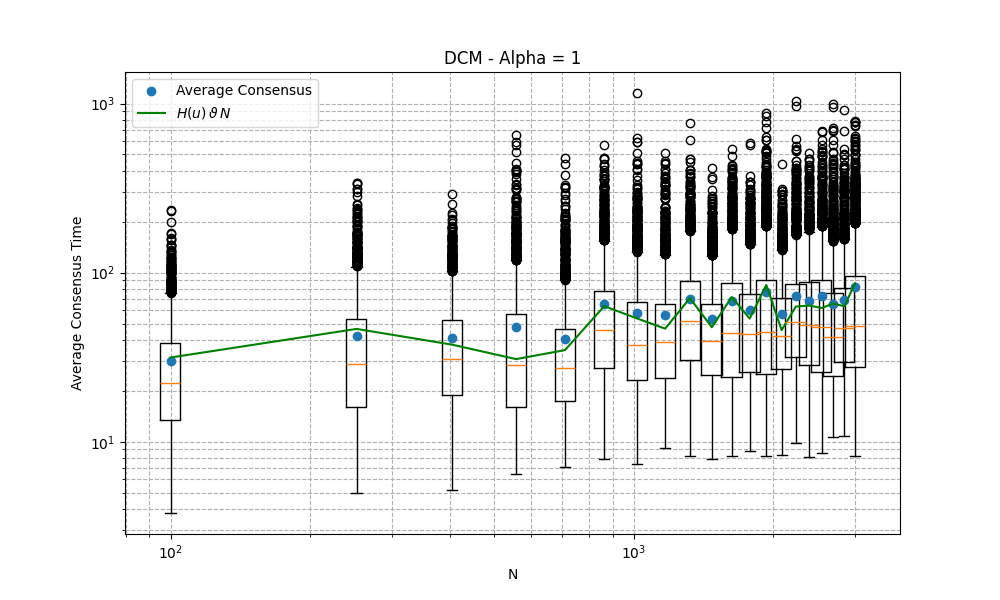}
        \caption{$\alpha = 1$}

    \end{subfigure}
    
    \vspace{0.5cm}  

    \begin{subfigure}{0.45\textwidth}
        \centering
        \includegraphics[width=\linewidth]{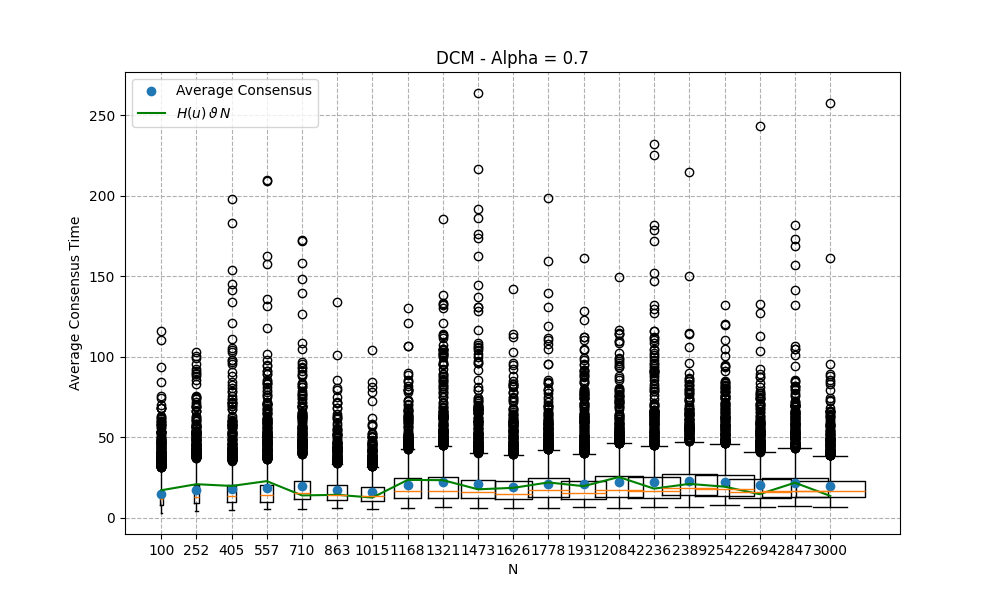}
        \caption{$\alpha = 0.7$}

    \end{subfigure}

    \caption{All figures represent the comparison between the box plots of the consensus time in the $\alpha$-DCM, for different ranges of $\alpha$ and increasing values of $n$. The blue dots denotes, for each $n$, the average over the voter model, the $\alpha$-DCM and the degree sequence realizations, while the box plots are taken over all the possible voter model outcomes. In figures (a)-(c) the plot is given in log-log scale and the green lines represent the predicted expected consensus time in terms of the conjectured value in \eqref{eq: conjecture genereal expected consensus}. In figures (d)-(e) the scale is not logarithmic as it pictures values that are slowly or not growing as $n$ increases. Here, for each $n$ ranging between 100 to 3000, we sampled 50 degree sequences, for each of which 20 graph realizations, for each of which 10 voter realizations.}
    \label{fig:DCM}
\end{figure*}
Consider $0<\alpha<\infty$, using an analogous argument as in the undirected case, we can explicitly compute the expressions in \eqref{eq:def-rho-gamma-delta-beta}, in terms of moments of the in- and out-degree sequences.  
Recall that, by construction, the sequences $d_x^\pm$ are supported in $[x_{\rm min}, \, d_{\rm max}^\pm]$, where $x_{\rm min}$ is the scale parameter of the Pareto random variable. Since we aim to impose strong connectivity with high probability \cite{CF04}, we choose $x_{\rm min} \geq 2$ both for the in- and out-degrees. Under the assumption for the $\alpha$-DCM, the empirical moments of the degree sequence $\Hat{m}^\pm_i$ can be approximated by their truncated moments
\begin{equation}
     m_i^\pm = \int_{x_{\rm min}}^{d^\pm_{\rm max}} x^i\,f_{D_\alpha}(x)\,dx\,.
\end{equation}
Observe that
\begin{equation} \label{eq:rho and gamma bounds}
    c_1\leq\rho\leq \frac{1}{x_{\rm min}}\,, \qquad c_2\leq\gamma\leq \frac{1}{x_{\rm min}}\beta\,,
\end{equation}
where
\begin{equation}
c_1 = \frac{x_{\rm min}}{\delta}\int_{x_{\rm min}}^{d^+_{\rm max}} \frac{1}{k} f_{D_\alpha}(k)\,dk = \Theta(\delta^{-1})\,, \ c_2 = x_{\rm min}\,c_1\,.
\end{equation}
As a consequence, we deduce the leading order of $\vartheta$ in \eqref{theta} to be
\begin{equation} \label{eq:delta over beta}
    \frac{\delta}{\beta} = \frac{(m_1^-)^2}{m_2^-}\,,
\end{equation}
for any $\alpha>0$. Note that the expression in \eqref{eq:delta over beta} involves only the first and second moments of the in-degree sequence $\mathbf{d^-}$, which is sampled according to $D_{\alpha}$, in the very same relation as the undirected case.  Indeed, if $\alpha >1$ then $\rho = \Theta(1)$ and bounded away from zero, while in the $\alpha\leq 1$ case $\rho$ can be vanishing, but we can easily expand 
$$\frac{1-\sqrt{1-\rho}}{\rho} = \frac{1}{2} + o(1)\,, \quad \text{as} \ \ \rho\to 0\,,$$
so that the first order asymptotic of the denominator of $\vartheta$ is never affected by the terms involving $\gamma$ and $\rho$. Therefore
\begin{equation} \label{eq:theta asymptotic}
    \vartheta = \Theta\Big((m_1^-)^2/m_2^-\Big)\,,
\end{equation}
as $n\to \infty$. This order is in accordance with the prediction in \cite{SAR08}. Moreover, the absence on the order of magnitude of the out-degree sequence can be interpreted by the fact that in the voter model on directed graphs the vertices adopt the opinion of out-neighbors, this forces the information to spread through the \emph{in-degrees}. On the other hand, the out-degrees become crucial in the expression of $\rho$ and $\gamma$ that determine the preconstant.

Finally, we analyze the components in \eqref{eq:def-rho-gamma-delta-beta} for different regimes of $\alpha>0$ using the bounds in \eqref{eq:rho and gamma bounds} and combine them in \eqref{theta}. Therefore, the predicted asymptotic order will be the same as those described in Subsection \ref{sec:refined results undirected}, which we report here for completeness.

\begin{table}[H]
    \centering
    \begin{tabular}{|c|c|c|c|c|c|}
        \hline
        & \textbf{$\alpha>2$} & \textbf{$\alpha =2$} & \textbf{$1<\alpha <2$}& \textbf{$\alpha =1$} & \textbf{$0<\alpha<1$} \\ 
        \hline
        $\delta$ & $\mathcal{O}(1)$ & $\mathcal{O}(1)$ & $\mathcal{O}(1)$ & $\log(n)$ & $n^\frac{1-\alpha}{\alpha}$\\ 
        \hline
        $\beta$ &  $\mathcal{O}(1)$ & $\log(n)$ & $n^\frac{2-\alpha}{\alpha}$ & $n$ & $n^\frac{2-\alpha}{\alpha}$\\ 
        \hline
        $\vartheta$ &  $\mathcal{O}(1)$ & $\log^{-1}(n)$ & $n^{-\frac{2-\alpha}{\alpha}}$ & $\log^2(n)\,n^{-1}$ & $n^{-1}$\\ 
        \hline
    \end{tabular}
    \label{tab:delta and beta}
\end{table}

Unlike in the configuration model setting, here we have the exact first-order behavior of $\mathbf{E}[\tau_{\rm cons}^u]$ for any regime $\alpha>0$. We have the precise preconstant that is encoded in the expression $\vartheta$ as $n$ tends to infinity. This is reflected in Figure \ref{fig:DCM}, where the accuracy of the predicted values and the empirical expected consensus remains relevant even for small values of $n$, and the expression is able to capture the fluctuation of the observable due to the variability of the graph structure and the degree sequence.\\
The box plot visualization in Figure \ref{fig:DCM} shows the presence of outliers, all
within the upper region above the third quartile. This evidence can be justified
by an intrinsic skewness of the consensus time distribution and some correlation
regarding the maximal degree of the related underlying network. The details of such observations are very similar to as described in Subsection \ref{sec:refined results undirected} for the undirected counterpart of the underlying geometry.


\subsubsection{Connectivity} So far we have always assumed the underlying (random) graph to be connected. In the case $\alpha$-CM and $\alpha$-DCM this assumption is translated into requiring the minimal (in-out) degree to be at least $2$ in the directed and $3$ in the undirected setting. We enforce it in the simulations by taking the minimal value of $\text{Pareto}(\alpha)$ large enough, i.e. $x_{\rm min}\geq 3$ (resp. $x_{\rm min}\geq 2$), so that the resulting random graph is connected with high probability \cite{FRdv17} (resp. strongly connected with high probability \cite{CF04}). Clearly the voter model is well defined on non-connected graphs, the only difference is that we need to redefine $\tau^u_{\rm cons}$ accordingly: it will be the first time such that all connected components reach local consensus. Therefore, on each connected component As a consequence, we conjecture that in such regimes a similar asymptotic order as in \eqref{eq:Sood prediction consensus} holds true, where $n$ is replaced by the size of the largest connected component (LCC), while the $m_i$ are replaced by the local moments of the degrees of the vertices inside the LCC. Similarly, in the directed framework we get an expression as \eqref{eq: conjecture genereal expected consensus} where $\vartheta$ and $n$ are restricted to the largest strongly connected component.

\section{Mean-field theory: Consensus \& Coalescence distribution}\label{sec: MFtheory: Coalescence}

The rigorous analysis of the voter model often rests on so‑called \emph{mean‐field conditions} for the underlying sequence of graphs. Informally, these conditions require that the random walk on $G_n$ to mix rapidly and that its stationary distribution is not too concentrated. We briefly recall the key definitions and limit theorems and then present our findings in the $\alpha$-CM and $\alpha$-DCM ensemble.

Thanks to the graphical construction mentioned earlier, the backward-in-time trajectories of the voter model can be directly linked to those of a system of coalescing random walks, where each lineage evolves according to the random walk generator defined in equation \eqref{CTRW}. Thus the characterization of $\tau^u_{\rm cons}$ turns to be essentially equivalent to that of the so-called \emph{coalescence time}, the first time the system of $n$ coalescing walks fully coalesce, that is:
\begin{equation}\label{1stmeetStat} \tau_{\rm coal}:= \inf\{t \geq 0 : X_t^{x}
= X_t^{y}, \, \forall x,y \in V\}, 
\end{equation}
where $(X_t^{x})_{t\geq 0}$ above is an independent copy of the Markov chain in \eqref{CTRW} starting from the vertex $x\in V$.
A crucial graph parameter to determine the law of $\tau_{\rm coal}$ is encoded in the \emph{mean meeting time for two random walks from stationarity}, i.e. : 
\begin{equation}
m_\pi:= \mathbf{E}[\tau^{\pi\times\pi}_{\rm meet} ]\,,
\end{equation}
where $\tau^{\pi\times\pi}_{\rm meet}:= \inf\{t \geq 0 : X_t^{x}
= X_t^{y}\, \mid x,y\overset{d}{=}\pi\}$ represents the first meeting time of two independent random walk lineages started according to the stationary distribution $\pi$. If the random walk \emph{mixing time} $t_{\rm mix}$ is small with respect to $m_\pi$, then any two lineages in the system of coalescing random walks will typically quickly equilibrate and their coalescence will be governed by the behaviour of the first meeting from stationarity.

Studying the distribution of the first meeting of two Markov chains from stationarity $\tau^{\pi\times\pi}_{\rm meet}$ can be translated into the question of studying the first hitting time of the diagonal set of the corresponding product chain. Asymptotically, the distribution of the latter scaled by its mean, turns to be exponentially distributed provided fast mixing of the chain (this corresponds to the classical Aldous-clumping heurstics, see \cite{AldPCH, MQS21}). This picture is valid for a sequence of complete graphs and has been made precise in other directed and undirected geometries in \cite{Oli12,Oli13}. In fact it was shown that asymptotically in $n$, the law of $\tau_{\rm coal}/m_\pi $ can be approximated by 
the law of an infinite sum  $\sum_{k> 1}Z_k$ of independent random variables $Z_k$'s with law \begin{equation}\label{eq:def-Z}
	Z_k\overset{d}{=}{\rm Exp}\left(\binom{k}{2}\right)\,,\qquad k\ge 2\, ,
\end{equation} provided either that the Markov chain in \eqref{CTRW} is \emph{transitive} and \emph{reversible} and such that 
\begin{equation}\label{MixBeforeMeet}
\lim_{n\to\infty}\frac{t_{\rm mix}}{m_\pi}=0,
\end{equation}
see \cite[Theorem 1.2]{Oli13}, or in the non-reversible setting, see \cite[Theorem 1.3]{Oli13}, provided  that 
 \begin{equation}\label{DirectedFastMixBeforeMeet}
\lim_{n\to\infty}(1 + q_{\rm max} t_{\rm mix})\pi_{\rm max}=0,
\end{equation}
with $q_{\rm max}:=\max_{x\in V} \sum_{y\neq x}q(x,y)$ and $\pi_{\rm max}=\max_{x\in V}\pi(x)$ being, respectively, the maxima of the exit rates from a state of an arbitrary generator as in \eqref{CTRW} and of the corresponding invariant measure $\pi$. Recall that we consider the voter model with $\eta_0\overset{d}{=} \text{Bern}(u)^{\otimes V}$, $u\in[0,1]$.
For the consensus time $\tau^u_{\rm cons}$, an analogous result holds true again either under \eqref{MixBeforeMeet} or under \eqref{DirectedFastMixBeforeMeet}(see \cite[Theorem 1.3]{Oli13}, and it states that asymptotically in $n$, the law of  $\tau^u_{\rm cons}/ m_\pi$ can be approximated by 
\begin{equation}\label{UE}
\sum_{k> K_u}Z_k
\end{equation}

where for $u\in (0,1)$
\begin{equation}\label{eq:def-K-oli}
\begin{split}
U&\overset{d}{=}{\rm Bern}(u)\,,\qquad A\overset{d}{=}{\rm Geom}(1-u)\,,\\
B&\overset{d}{=}{\rm Geom}(u)\,, \qquad K_u\overset{d}{=} U A + (1-U) B\,,
\end{split}
\end{equation}
and the $Z_k$'s are independent random variables with law as in \eqref{eq:def-Z}.
The fact that the sum starts from $K_u$ rather than $2$, as for $\tau_{\rm coal}$,
captures the initialization of the voter opinions for which achieving consensus can be driven by $K_u$ ancestors rather than the last common ancestor of the backward genealogy. The convergence in \eqref{UE} was provided in Wasserstein metric on probability measures, and this guaranteed the convergence of the corresponding means, and hence that under the above mean-field conditions, the expected consensus grows twice as the meeting:
\begin{equation}\label{2Meeting}
\E[\tau^u_{\rm cons}]\sim 2 H(u) m_\pi.\end{equation}

\begin{figure*}[t] 
    \centering
    \begin{subfigure}{0.45\textwidth}
        \centering
        \includegraphics[width=\linewidth]{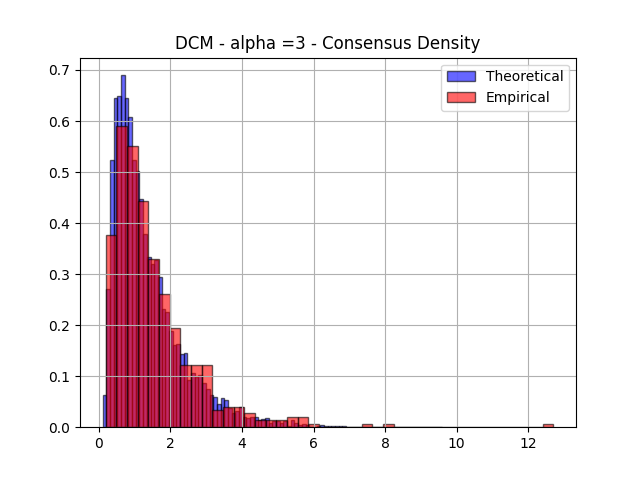}
        \caption{$\alpha = 3$ }

    \end{subfigure}
    \hfill
    \begin{subfigure}{0.45\textwidth}
        \centering
        \includegraphics[width=\linewidth]{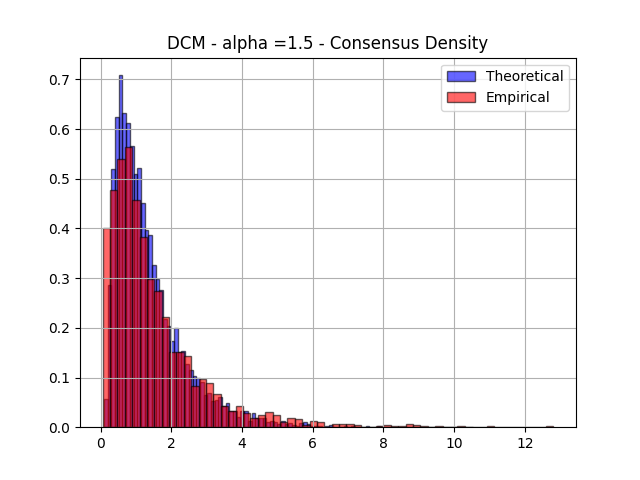}
        \caption{$\alpha = 1.5$}

    \end{subfigure}
    
    \vspace{0.5cm}  

    \begin{subfigure}{0.45\textwidth}
        \centering
        \includegraphics[width=\linewidth]{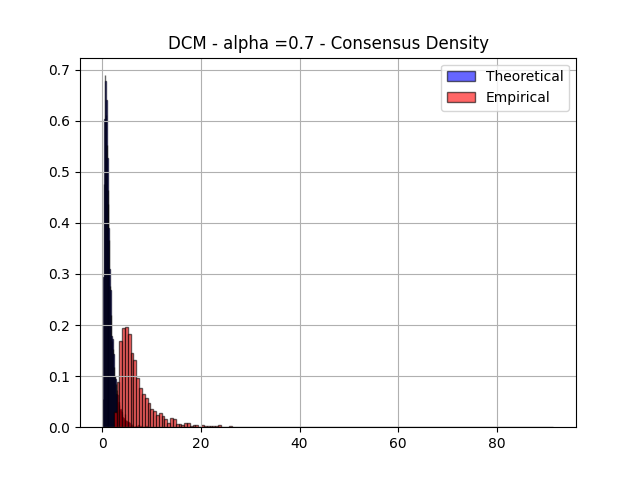}
        \caption{$\alpha = 0.7$}

    \end{subfigure}

    \caption{In blue, the numerical expression for the target density of sum of exponential random variables described in \eqref{UE}. In red, the numerical density of the rescaled consensus time in the  $\alpha$-DCM ensemble for different values of $\alpha$. We considered graphs with $n=1000$ vertices, and run 50 voter model iterations for 5 different graphs realizations.}
    \label{fig:Density_consensus_DCM}
\end{figure*}

In practical terms, this implies that, when mean-field conditions hold, determining the asymptotic behavior of $\E[\tau^u_{\rm cons}]$ equals (up to a factor of two) to analyzing the asymptotic behavior of $m_\pi$. This relationship forms the basis for the conjecture in \eqref{eq: conjecture genereal expected consensus}, which was explicitly verified for the directed graph model with deterministic bounded degrees in \cite{ACHQ23}. In this directed ensemble, the condition in \eqref{DirectedFastMixBeforeMeet} can be checked to hold true with high probability despite there is no closed formula for the invariant measure $\pi$. Moreover, since such an ensemble is locally tree like up to the $\log n$ scale of the mixing, one can work out the underlying combinatorics via path counting estimates on directed trees and obtain the explicit formula in \eqref{theta}. In particular, we stress that, in light of \eqref{2Meeting}, our conjecture in \eqref{eq: conjecture genereal expected consensus} for the $\alpha$-DCM ensemble, suggests that when the condition in \eqref{DirectedFastMixBeforeMeet} is satisfied, it is natural to predict that with high probability w.r.t. $\mathbb{P}$ 
\begin{equation}\label{MeetingConjecture}
    m_\pi\sim n\vartheta /2 .  
\end{equation}
 
In the undirected analogous locally-tree-like bounded degree case, although the invariant measure is explicit and the reversible mean-field condition in \eqref{MixBeforeMeet} can be verified, the combinatorics related to the random walk Green's function becomes structurally different. This increased complexity, compared to the directed case, arises from the numerous backtrackings of the random walk.
 As a result, except for $d$–regular undirected cases such as \cite{CF04} no formula has been found out, predicting the precise first-order remains a difficult open problem. We now explore via simulations whether the coalescence result state above in \eqref{UE} is valid in the directed case, specifically the $\alpha-$DCM ensemble. 

\paragraph{\bf Coalescence in directed Pareto ensemble with finite mean degrees}\label{Coalescence}

In Figure~\ref{fig:Density_consensus_DCM}, we present simulations comparing the target density, derived from the sum of independent exponential random variables suggested by the mean-field approximation in equation~\eqref{UE}, with the empirical density of the rescaled consensus time across different regimes of the $\alpha$-DCM ensemble. 
Our findings indicate that when the degree distribution has a finite mean (i.e., $\alpha > 1$), the empirical and theoretical densities closely agree, with any discrepancy being negligible in the large-size limit.
In contrast, when the degree distribution has infinite mean (i.e., $\alpha \leq 1$), the shape of the empirical density deviates noticeably from the target distribution, although the two remain qualitatively similar. This suggests that while the mean-field approximation captures essential features of the consensus dynamics, it does not fully account for the additional variability introduced by heavy-tailed degrees.

In conclusion, our simulations strongly support mean-field approximations for $\alpha > 1$ but also clearly illustrate their limitations when $\alpha \leq 1$. This highlights the need for refined theoretical models to accommodate the additional variability and correlation structure observed in the infinite-mean setting.

\section{Mean-field theory: Wright-Fisher diffusive approximation}\label{sec: MFtheory: diffusions}

\begin{figure*}[t]
    \centering
    \begin{subfigure}{0.45\textwidth}
        \centering
        \includegraphics[width=\linewidth]{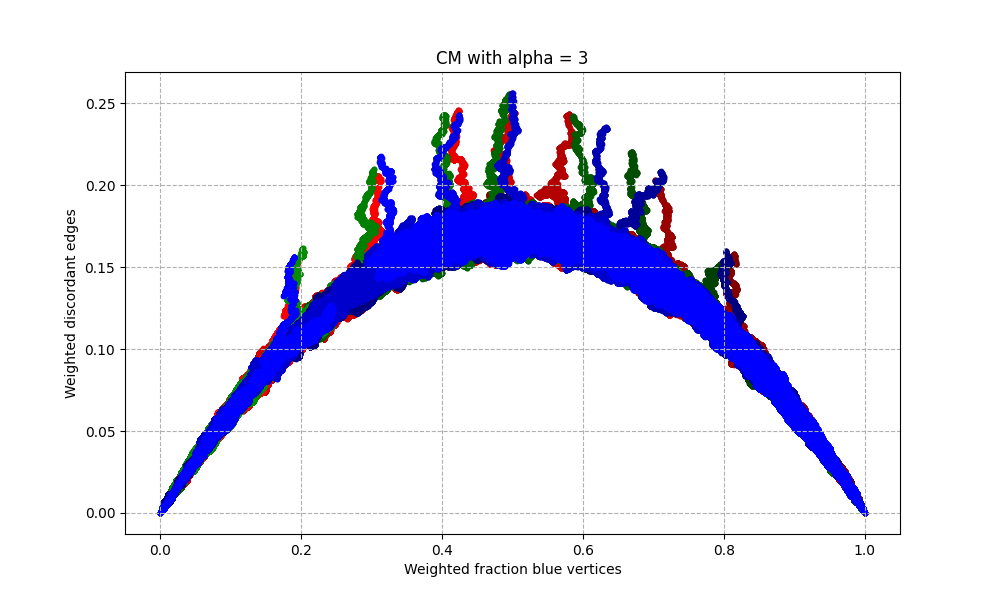}
        \caption{$\alpha=3$}
  
    \end{subfigure}
    \hfill
    \begin{subfigure}{0.45\textwidth}
        \centering
        \includegraphics[width=\linewidth]{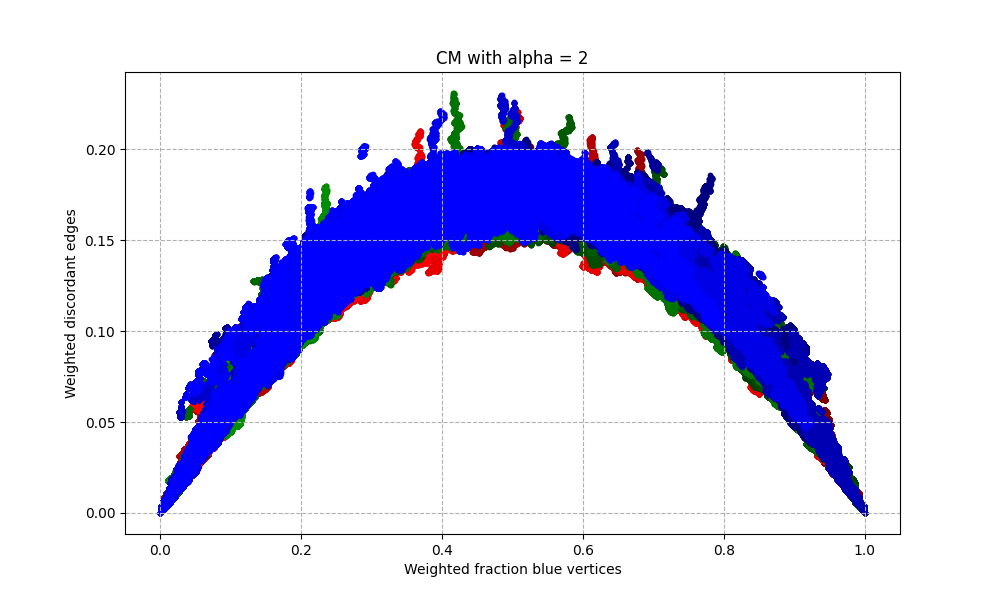}
        \caption{$\alpha=2$}
  
    \end{subfigure}
    

    \begin{subfigure}{0.45\textwidth}
        \centering
        \includegraphics[width=\linewidth]{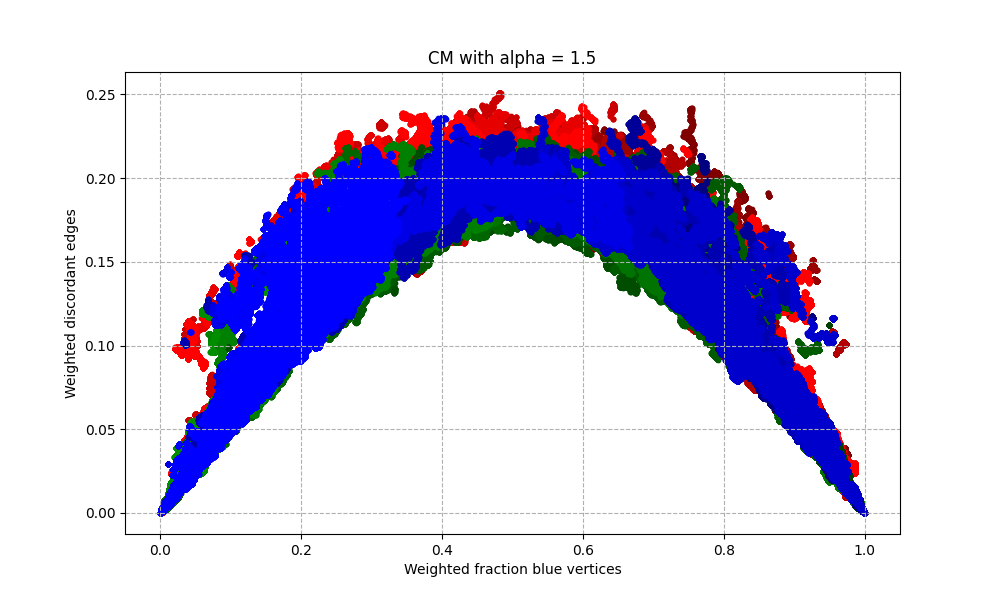}
        \caption{$\alpha=1.5$}
  
    \end{subfigure}
    \hfill
    \begin{subfigure}{0.45\textwidth}
        \centering
        \includegraphics[width=\linewidth]{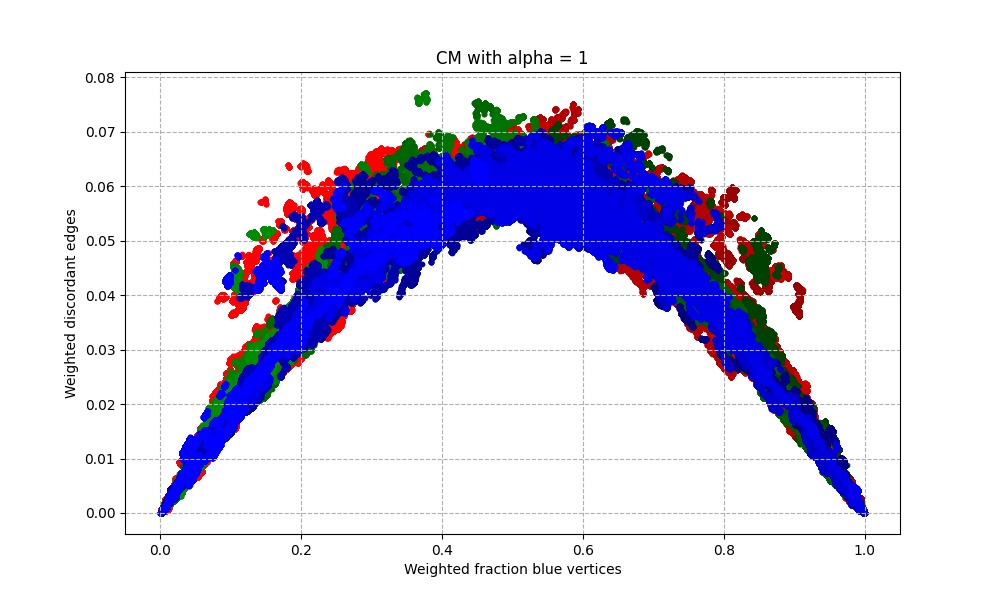}
        \caption{$\alpha=1$}
  
    \end{subfigure}
        

    \begin{subfigure}{0.45\textwidth}
        \centering
        \includegraphics[width=\linewidth]{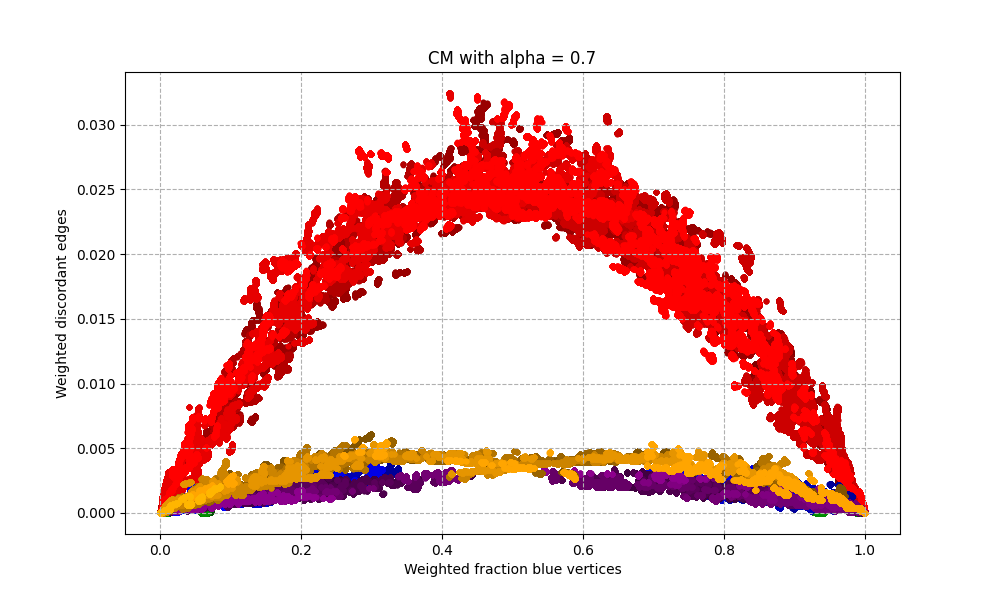}
        \caption{$\alpha=0.7$}
  
    \end{subfigure}

    \caption{For each figure, on the $x$ axes the $\pi$-weighted density of opinion 1, while on the $y$ axes the predictable quadratic variation $\langle M_n\rangle_t$. For figures (a)-(d), the different colors represent different degree sequences and corresponding graphs sampled according $\alpha$-CM, while different shades of the same color represent different initialization parameter $u$ for the product of Bernoulli random variables. For figure (e) the degree sequence is quenched and different colors represent different graph realizations.} 
    \label{fig:WF-Undirected}
\end{figure*}
The mean-field conditions discussed in the previous section are not only useful to
describe the ancestry evolution of the opinions but also to give a diffusive approximation for the evolution of the density of opinions of one type which we next briefly describe and again test in the $\alpha$-CM and $\alpha$-DCM case. 

Let us first consider the \emph{empirical density} of the opinions of type $1$:
\begin{equation}\label{OpDensity}
  \mathcal{O}_n(t)=\frac1n\sum_{x\in V}\eta_t(x).
\end{equation} 

On the complete graph it is easy to check that $(\mathcal{O}_n(t))_{t\geq 0}$ is a mean-zero martingale with jumps of size $1/n$, whose predictable quadratic variation is of the form $\frac{n}{n-1}\mathcal{O}_n(t)(1-\mathcal{O}_n(t))$. We recall that if $(M_t)_{t\ge 0}$ is a square integrable continuous time martingale then its predictable quadratic variation, denoted as $\langle M\rangle_t$, is the unique predictable, increasing process such that $M_t^2-\langle M\rangle_t$ is also a martingale. This process tracks the expected accumulated variance based only on past information. When sped-up by a factor $n$  (which is the leading order of the consensus and of $2m_\pi$ in this setting), the trajectories of the martingale $\mathcal{O}_n(tn)$ converges to the so-called standard \emph{Wright-Fisher} diffusion $(Y_t)_{t\ge 0}$. Such a diffusion is a fundamental stochastic process in population genetics that models the evolution of allele frequencies under genetic drift, and it is characterized as the unique solution to the martingale problem associated with the generator $\frac{x(1 - x)}{2} \frac{d^2}{dx^2}$. Equivalently, $(Y_t)_{t\geq 0}$ can be described as the unique diffusion process on $[0,1]$ satisfying the SDE
\begin{equation}\label{WF}
dY_t = \sqrt{Y_t(1 - Y_t)} dW_t,
\end{equation}
where $ (W_t)_{t\geq 0}$ is a standard Brownian motion from any starting value $W_0=u\in[0,1]$.

In particular, if a square-integrable discrete martingale sequence $M_n(t_n)$ with values in $[0,1]$ and vanishing increments has predictable quadratic variation of the form $\int_0^{t_n} M_n(s)(1-M_n(s)ds$, then the Wright-Fisher approximation kicks in.

The key works in \cite{CCC16,Che18} study whether the same diffusive approximation hold true for the classical voter model or more general variants under mean field conditions analogous to those discussed in the previous section. For the standard voter model in \eqref{sec:VoterDef}, as soon as we move out of the $d$-regular undirected setting, it is no longer true in general that the density of opinion in \eqref{OpDensity} is a martingale. Yet, it is still the case (see \cite{CCC16}) that in an irreducible setting for which there is a unique equilibrium $\pi$ associated with \eqref{CTRW}, the $\pi$-\emph{weighted density} of the opinions of type $1$ 
\begin{equation} \label{eq:barpn}
M_n(t):=\sum_{x\in V}\pi(x)\,\eta_t(x)\,,
\end{equation}
is again a martingale with jump sizes bounded by $\pi_{\rm max}:=\displaystyle\max_{x\in V}\pi(x)$ and whose predictable quadratic variation $\bigl\langle M_n \bigr\rangle_t$ can be expressed as
\begin{equation}\label{PQV}
\begin{split}
\frac{d}{dt} \bigl\langle M_n \bigr\rangle_t&= \sum_{\substack{x, y\in V \\ x\neq y}} \pi^2(x) q(x, y)\left[\eta_t(x)+\eta_t(y)-2 \eta_t(x) \eta_t(y)\right]\\
&=\pi_{\Delta}\sum_{x\in V}\frac{\pi(x)^2}{\pi_{\Delta}}\frac{\vec{D}_x(t)}{d^+_x}\,,
\end{split}
\end{equation}
where in the second expression, valid for the RW kernel in \eqref{eq:VM rates},
$\vec{D}_x(t)$ stands for the number of edges outgoing vertex $x$ with different opinions on their extremes and $\pi_{\Delta}:=\sum_{x\in V}\pi(x)^2$.
This second expression says that the quadratic variation grows as 
$\pi_\Delta$ and it is captured by the local densities $\vec{D}_x(t)/d^+_x$ of the \emph{discordant edges} weighted with the normalized \emph{squared invariant masses}  $\pi(x)^2/\pi_{\Delta}$ (see e.g. \cite{ABHHQ22,avena2025, FdHreview} for more details).  
The main result in \cite{CCC16}, see Theorem 2.2 therein, states that 
for any RW kernel $q$, if the invariant is not too concentrated in the sense that
\begin{equation}\label{PiDiag}\lim_{n\to\infty}\pi_{\Delta}=0,
\end{equation}
 and if the fast mixing condition in \eqref{MixBeforeMeet} is valid, then, as for the complete graph, the weighted-density process $M_n(t\times m_\pi)$ sped-up by the mean meeting time sequence $m_\pi$ converges to the Wright-Fisher diffusion in \eqref{WF}.
In particular the latter is true since for any $t>0$, on time scale $m_\pi$, the predictable quadratic variation in \eqref{PQV} behaves as discussed right after \eqref{WF} (see also \cite[Theorem 2.1]{CCC16}), that is:

\begin{equation}\label{Chen}
\begin{split}
\bigl\langle M_n \bigr\rangle_{t\times m_\pi}&=m_\pi\pi_{\Delta}\int_0^{t}\sum_{x\in V}\frac{\pi(x)^2}{\pi_{\Delta}}\frac{\vec{D}_x(s)}{d^+_x}ds \\
&\sim \int_0^{t} M_n(s\times m_\pi)\left(1-M_n(s\times m_\pi)\right)ds.    
\end{split}
\end{equation}

To see why the latter is true, a rough heuristic is as follows. For fast-mixing graphs in the sense of condition \eqref{MixBeforeMeet} at the consensus time scale $t_n$ local correlations are weak, and we may guess that
\begin{equation*}
\begin{split}
     \eta_{t_n}(x)&\sim \mathbf{E}[\eta_{t_n}(x)]  \quad \text{and} \\
\mathbf{E}[\eta_{t_n}(x) \eta_{t_n}(y)] &\sim \mathbf{E}[\eta_{t_n}(x)]\mathbf{E}[\eta_{t_n}(y)].
\end{split}
\end{equation*}
 Plugging such  estimates in the sum in the first expression in \eqref{PQV}, we see that  
\begin{equation}\notag
  d\bigl\langle M_n \bigr\rangle_{t_n}
\approx M_n(t_n)\bigl(1-M_n(t_n)\bigr) dt,
\end{equation} as in \eqref{Chen}.

\section{ W-F in Pareto ensembles with finite mean \& effective diffusion parameter}

\begin{figure*}[t]
    \centering
    \begin{subfigure}{0.45\textwidth}
        \centering
        \includegraphics[width=\linewidth]{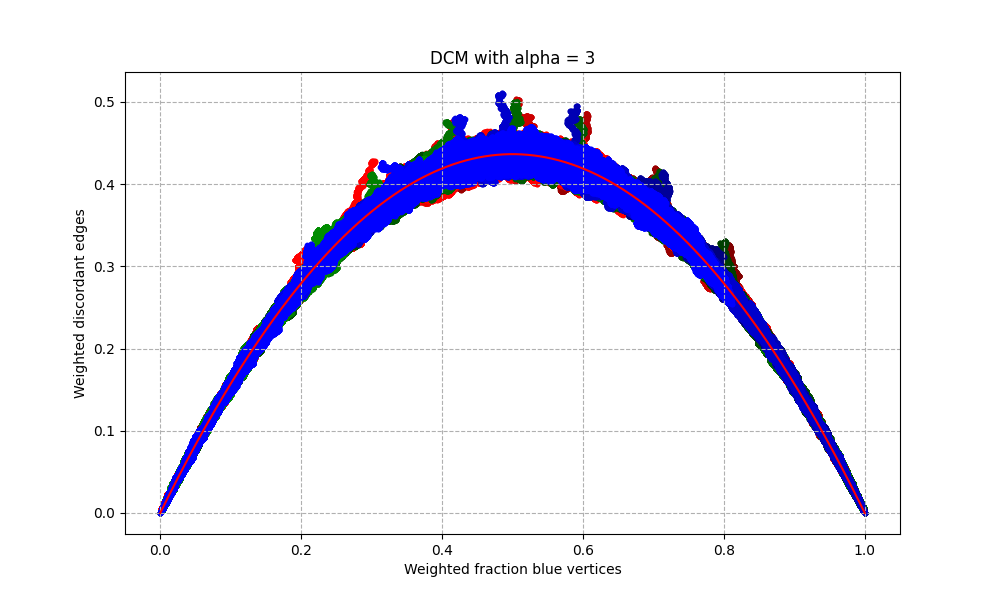}
        \caption{$\alpha=3$}
  
    \end{subfigure}
    \hfill
    \begin{subfigure}{0.45\textwidth}
        \centering
        \includegraphics[width=\linewidth]{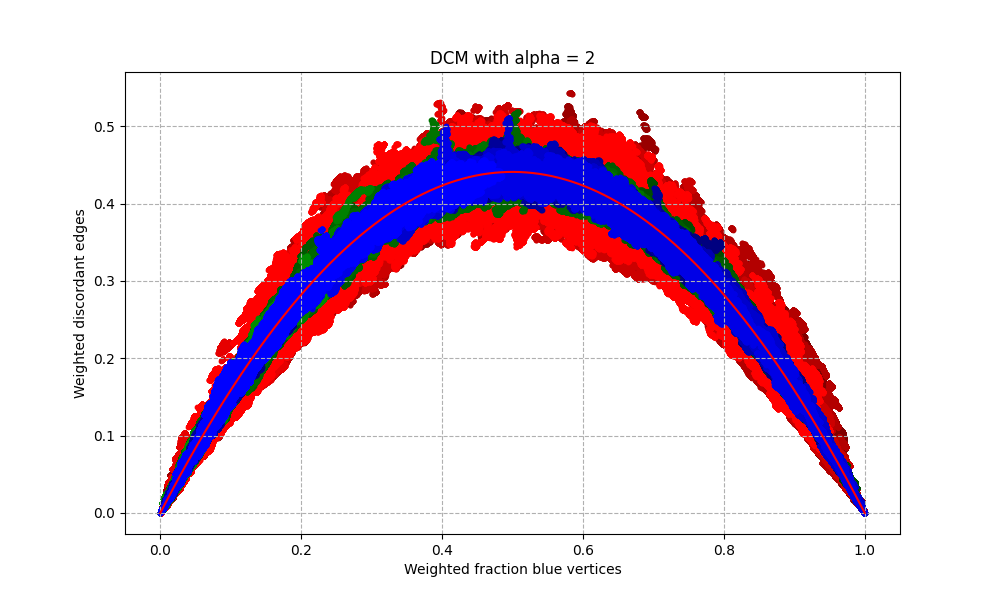}
        \caption{$\alpha=2$}
  
    \end{subfigure}
    

    \begin{subfigure}{0.45\textwidth}
        \centering
        \includegraphics[width=\linewidth]{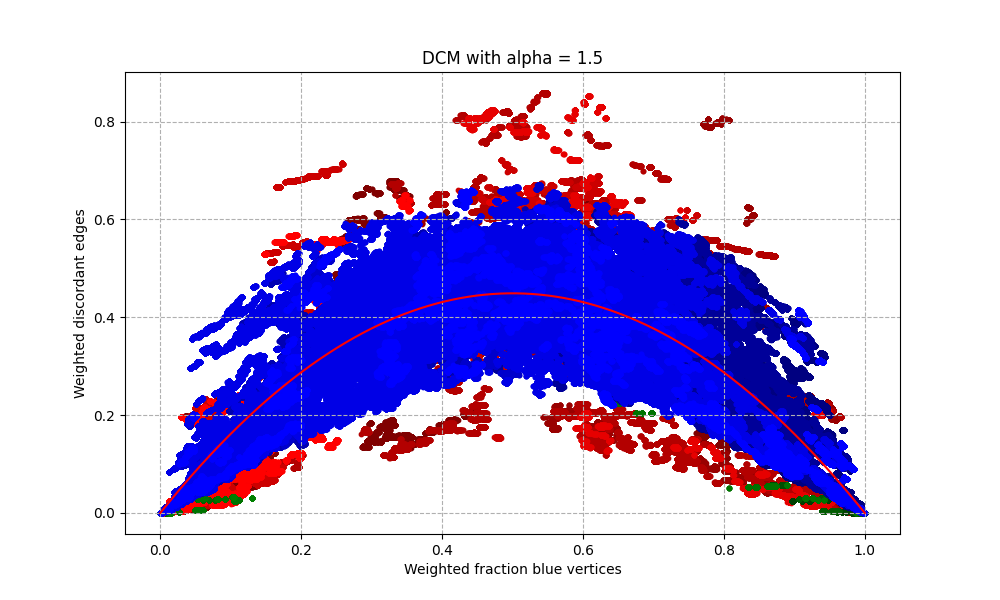}
        \caption{$\alpha=1.5$}
  
    \end{subfigure}
    \hfill
    \begin{subfigure}{0.45\textwidth}
        \centering
        \includegraphics[width=\linewidth]{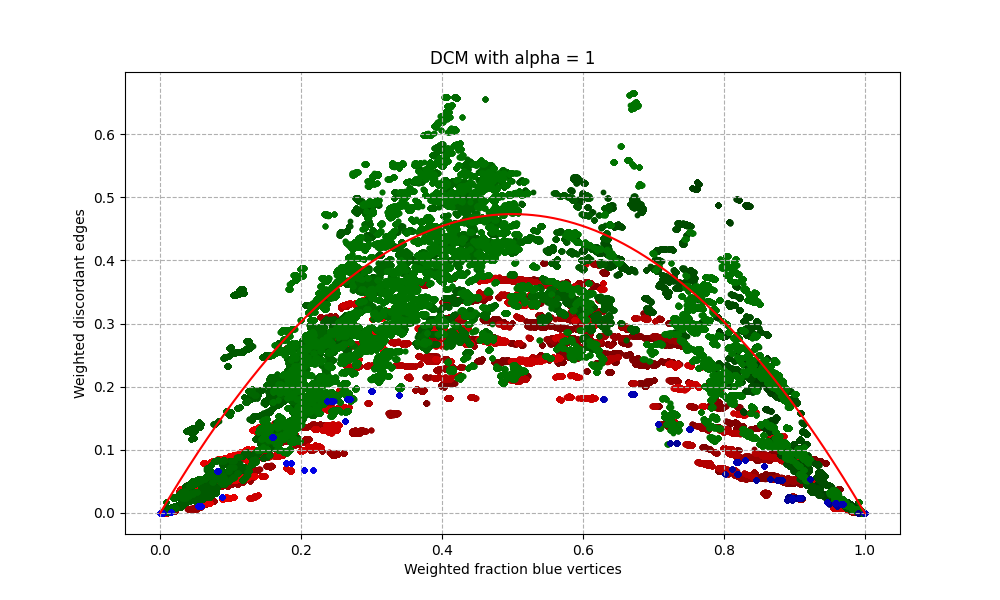}
        \caption{$\alpha=1$}
  
    \end{subfigure}
        

    \begin{subfigure}{0.45\textwidth}
        \centering
        \includegraphics[width=\linewidth]{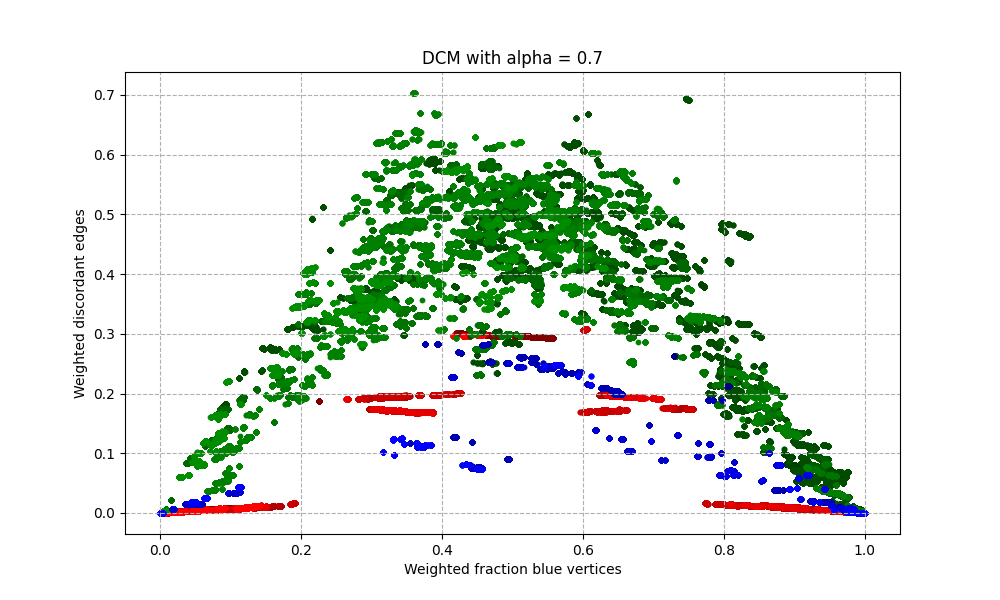}
        \caption{$\alpha=0.7$}
  
    \end{subfigure}

    \caption{For each figure, on the $x$ axes the $\pi$-weighted density of opinion 1, while on the $y$ axes the predictable quadratic variation $\langle M_n\rangle_t$. The different colors represent different degree sequences and corresponding graphs sampled according $\alpha$-DCM, while different shades of the same color represent different initialization parameter $u$ for the product of Bernoulli random variables. The red parabola in (a)-(d) is $y=\chi\,x(1-x)$, where $\chi$ is the diffusion parameter in \eqref{eq: chi DCM}}.
    \label{fig:WF-Directed}
\end{figure*}

We validate the Wright–Fisher approximation with simulations on both the $\alpha$-CM and $\alpha$-DCM models with $\alpha\in \{ 3, 2, 1.5, 1, 0.7\}$, respectively. The findings for the undirected case are reported in  Figure \ref{fig:WF-Undirected}, where different colors represent different realizations of the degree sequences and corresponding graphs, while different shades of the same color represent different
initialization parameter $u$ for the product of Bernoulli random variables. 

In these simulations, we plot the 
weighted opinion density $M_n(t)$ on the $x$-axis and the \emph{normalized weighted discordances } $\displaystyle\sum_{x\in V}\frac{\pi(x)^2}{\pi_{\Delta}} \frac{\vec{D}_x(t)}{d^+_x}$ on the $y$-axis, or equivalently, the derivative of the predictable quadratic variation in \eqref{Chen} up to the meeting time factor $m_{\pi}$. 
The resulting curve clearly exhibits parabolic behavior, passes through 0 and 1, and peaks when the weighted opinion density is near $0.5$ in all cases where $\alpha>1$. That is, in all these tail-parametric regimes our simulation data closely follow a parabolic curve, as shown by the representative points lining up along $y=\chi x(1-x)$ for some $\chi \in (0,1]$.
This \emph{effective diffusion parameter} $\chi$ emerging in the limit can be interpreted in terms of structural properties of the network closely tied to the meeting time $m_\pi$ and the stationary mass concentration $\pi_\Delta$ and it is in fact captured by the relation
\begin{equation} \label{eq:chi}
\frac{1}{\chi}\sim m_\pi\pi_{\Delta},
\end{equation}
encoded in the first expression in \eqref{Chen}.

Heuristically, a large $m_\pi$ (slow random-walk meeting) or a large $\pi_\Delta$ (highly uneven stationary distribution) will reduce the effective variance yielding a smaller $\chi$. 
For example in a complete graph, $\pi(x)=1/n$ and it is known that $m_\pi\sim n$, giving $m_\pi \pi_\Delta\approx 1$ and thus $\chi=1$. In contrast, for a heterogeneous network, $\pi_\Delta$ can be significantly larger (if a few high-degree nodes carry a large stationary mass), and $m_\pi$ may also scale differently. This typically leads to $\chi<1$, meaning the diffusion has lower variance than the dense well-mixed case.

Empirically, as long as the degree distribution has a finite mean and  $\pi_{\Delta}\to 0$ the voter model’s behavior aligns well with the Wright–Fisher diffusion. The macroscopic opinion $M_n(t)$ diffuses gradually, and the plot of the normalized quadratic variation remains close to the theoretical parabola. As evident from Figure \ref{fig:WF-Undirected}, deviations emerge only in regimes where the mean-field conditions break down in extremely heterogeneous networks with infinite expected degree (such as power-law with exponent $\alpha\le 1$), where one or a few nodes hold a non-negligible fraction of $\pi$ and drive the global consensus formation (recall that in a star graph it can be checked that the WF approximation does not hold). 

The directed case is completely analogous, as reported in Figure \ref{fig:WF-Directed}, where a similar trend can be observed, although the curves show more variability. This is expected due to the randomness of $\pi$ in the case of DCM. Still for $\alpha>1$, the quadratic variation tracks the parabola, confirming the presence of Wright-Fisher diffusion.  Deviations are more pronounced than in the undirected case, due to the added randomness and lack of symmetry.

Furthermore, in the DCM case we manage to express the explicit value of $\chi$, despite the randomness of $\pi_\Delta$, in terms of the in- and out-degree sequences. In fact, we get
\begin{equation} \label{eq: chi DCM}
    \chi = 1-\frac{1-\sqrt{1-\rho}}{\delta\,\rho}\,,
\end{equation}
where $\delta,\rho$ as in \eqref{eq:def-rho-gamma-delta-beta}. This value comes as a consequence of the expression in \eqref{eq:chi} and some computations for the expected density of discordant edges in \cite{Cap25} under the bounded degree assumption. As a consequence, we can interpret the diffusion parameter $\chi$ as a metastable state at which the density of discordant edges stabilizes for a long time before reaching consensus. We can now compare our prediction in \eqref{eq: conjecture genereal expected consensus} for the expected consensus time in general DCM models with the latter expression for $\chi$, and conjecture that the effective diffusion parameter $\chi$ has the form expressed in \eqref{eq: chi DCM} for any DCM realization with finite-mean degree distribution. We validate this claim via simulations on the $\alpha$-DCM ensemble for different values of $\alpha$ shown in Figure \ref{fig:WF-Directed}, and notice that as long as $\alpha>1$ the red parabola $y=\chi\,x(1-x)$ matches the prediction. 

In conclusion, both in the directed and the undirected ensembles, apart from the extreme degree volatile cases corresponding to $\alpha<1$, the Wright–Fisher diffusion approximation provides a remarkably accurate and analytically tractable description of the voter model’s evolution on heterogeneous networks on the consensus time-scale (recall \eqref{2Meeting}). It connects the micro-level events (through discordant edges, see $\vec{D}_x(t)$ in \eqref{Chen} ) to a macro-level diffusion law, and it highlights how network topology enters via the effective variance parameter~$\chi$. 


\bibliography{main.bib}

\end{document}